\documentclass[preprint,12pt]{elsarticle}
\usepackage{amsmath,amssymb,amsfonts,bm}
\usepackage{graphicx,booktabs,multirow}
\usepackage{xcolor,hyperref}
\usepackage[capitalize,noabbrev]{cleveref}
\usepackage{enumitem}
\usepackage{acro}
\usepackage{algorithm}
\usepackage{algpseudocode}
\usepackage{amsmath}

\graphicspath{{../numerical_examples/}}

\usepackage{amsthm}
\theoremstyle{remark}
\newtheorem{remark}{Remark}[section]
\newtheorem{property}{Property}[section]
\usepackage{comment}
\usepackage{multirow}

\acsetup{make-links = false}

\DeclareAcronym{KED}{
  short = KED,
  long  = kinetic defect measure
}

\DeclareAcronym{HCL}{
  short = HCL,
  long  = hyperbolic conservation laws
}

\DeclareAcronym{ROM}{
  short = ROM,
  long  = reduced order model,
  short-plural-form = ROMs,
  long-plural-form  = reduced order models
}

\DeclareAcronym{DMD}{
  short = DMD,
  long  = dynamic mode decomposition  
}

\begin{document}
\begin{frontmatter}
\title{Structure-Informed Data-Driven Reduced-Order Modeling of Scalar Hyperbolic Conservation Laws via Kinetic Defect Measure}
\author[oden]{Marissa Llamas}
\author[oden,ase]{Jan Fuhg}
\author[oden,ase]{Hannah Lu\corref{cor1}}

\ead{hannah.lu@austin.utexas.edu}

\cortext[cor1]{Corresponding author}

\address[oden]{
Oden Institute for Computational Engineering and Sciences,
The University of Texas at Austin,
201 E 24th St, Austin, TX 78712, USA
}

\address[ase]{
Department of Aerospace Engineering and Engineering Mechanics,
The University of Texas at Austin,
2617 Wichita Street, Austin, TX 78712, USA
}

\begin{abstract}
Reduced-order modeling of transport-dominated systems remains challenging because moving fronts and shocks are poorly represented by low-dimensional linear subspaces. We develop a structure-informed data-driven \ac{ROM} for scalar hyperbolic conservation laws based on the kinetic defect formulation. This formulation separates the nonlinear dynamics into known characteristic transport and a kinetic entropy defect localized on the shock manifold. We exploit this structure by first removing the known transport from the solution snapshots. We then extract and register the remaining defect-driven dynamics in a shock-attached coordinate system. Separate \acp{ROM} are used to evolve the shock geometry and the registered defect-driven source. During prediction, the predicted shock geometry is used to inverse-register the learned defect-driven source, which advances the kinetic state and recovers the physical solution. Numerical examples in one and two spatial dimensions demonstrate accurate reconstruction and prediction of nonlinear transport with shocks, including evolution beyond the training interval, while accurately capturing the mass and entropy-dissipation behavior of the reference solution.
\end{abstract}

\begin{keyword}
Reduced-order Models \sep Dynamic Mode Decomposition \sep Kinetic Formulation \sep Entropy Defect  \sep Transport-dominant Systems\sep Conservation Laws
\end{keyword}
\end{frontmatter}

\section{Introduction}\label{sec:intro}
Reduced order models (ROMs) provide computationally efficient representations of high-dimensional dynamical systems and are particularly valuable in repeated-query applications such as uncertainty quantification, optimization, data assimilation, and control. Classical approaches, including proper orthogonal decomposition (POD)~\cite{benner2015survey,berkooz1993proper} and \ac{DMD}~\cite{schmid2010dynamic,kutz2016dynamic}, seek low-dimensional representations of the dominant solution dynamics from high-fidelity snapshots. Their application to transport-dominated problems, however, remains challenging. Solution manifolds associated with moving fronts and shocks often exhibit slowly decaying Kolmogorov $n$-widths~\cite{ohlberger2016reduced,greif2019decay, peherstorfer2022breaking,rim2023manifold}: even when the underlying dynamics are simple, translations of a localized structure generally require many spatial modes. Consequently, conventional Eulerian ROMs may require large reduced dimensions to represent transport and can produce substantial errors or spurious oscillations under aggressive truncation~\cite{lu2020lagrangian}. This difficulty is compounded for nonlinear hyperbolic conservation laws, where intersecting characteristics generate shocks, and the reduced model must represent not only transport but also shock formation and propagation while recovering the physically admissible entropy solution.

Several strategies have been developed to overcome this limitation by seeking representations in which transport-dominated dynamics are more compressible. Transport-aware approaches based on local or adaptive bases, solution alignment, registration, and nonlinear coordinate transformations can remove or reduce the apparent motion of coherent structures before reduction~\cite{ohlberger2013nonlinear,reiss2018shifted,iollo2014advection,rim2018displacement,taddei2020registration,rim2023manifold,peherstorfer2015online,peherstorfer2020model}, while nonlinear-manifold approaches~\cite{lee2020model,ehrlacher2020nonlinear,taddei2015reduced,rim2018transport,barnett2022quadratic,barnett2023neural,kim2022fast}, including autoencoder-based ROMs, replace the linear approximation space with a nonlinear approximation manifold. Such methods can represent moving and discontinuous structures more efficiently than global linear subspaces, but generally require sufficiently rich training data and/or introduce additional complexity in constructing and evolving the reduced representation. An alternative is to exploit known transport structure through a change of coordinates. Lagrangian POD and DMD follow characteristic trajectories and can recover low-rank representations for smooth advection-dominated problems \cite{mojgani2017lagrangian,lu2020lagrangian}. Once shocks form, however, characteristics intersect, and the Lagrangian mapping ceases to be one-to-one. The resulting grid distortion can become particularly severe after projection onto a reduced space, where the topology of the intersecting characteristic map is not necessarily preserved \cite{lu2021dynamic}. Alternative physics-aware transformations can incorporate shock information explicitly. In particular, hodograph-based \ac{DMD}~\cite{lu2021dynamic} can accommodate shocks by interchanging dependent and independent variables and explicitly incorporating shock information, but its construction relies on invertibility or decomposition into monotone solution branches and does not extend straightforwardly to multiple spatial dimensions. These limitations motivate a representation that retains the efficiency of physics-aware reduction while treating shock dynamics in a more intrinsic manner.

The kinetic formulation of scalar hyperbolic conservation laws provides such a representation~\cite{perthame1991kinetic, lions1994kinetic,perthame2002kinetic}. By introducing an auxiliary kinetic coordinate, the \emph{nonlinear} conservation law is lifted to a higher-dimensional \emph{linear} transport equation with known characteristic velocity. The departure from this free transport is represented by a kinetic entropy defect measure, which accounts for entropy production associated with shocks~\cite{lions1994kinetic,perthame2002kinetic}. The physical solution is recovered exactly by integration over the kinetic coordinate. The kinetic formulation therefore provides a natural decomposition of the nonlinear dynamics into a known transport component and a defect-driven correction.

This decomposition is particularly attractive for model reduction because the kinetic defect possesses strong structural properties~\cite{perthame1998uniqueness,perthame2002kinetic}. It is nonnegative, has bounded support in the kinetic coordinate, and vanishes wherever the entropy solution is smooth. For piecewise smooth solutions, its nontrivial spatial support is confined to the lower-dimensional shock manifold. Recent work has shown that this localization can be exploited to represent and evaluate the defect through quantities defined along shock trajectories, including in multiple spatial dimensions~\cite{srivastava2026computable}. This localization suggests that model reduction can focus on the dynamically active departure from free transport rather than on the complete translating solution field.

Motivated by this observation, we develop a \emph{structure-informed data-driven \ac{ROM} for scalar hyperbolic conservation laws via the kinetic defect measure}. Solution snapshots are first lifted to kinetic space and pulled back along the analytically known characteristics, removing the transport prescribed by the flux. Temporal changes in this characteristic frame provide an empirical approximation of the defect-driven evolution. We then exploit its localization by extracting the evolving shock manifold and registering the defect-driven evolution in a shock-attached coordinate system. The \acp{ROM} for shock geometry and the registered defect-driven source are learned separately. During online prediction, the predicted shock geometry is used to inverse-register the learned defect-driven source, which is combined with the prescribed characteristic transport to advance the kinetic solution and recover the physical state. In this way, the governing structure determines \emph{what} must be learned: the known transport is prescribed, while model reduction is restricted to the localized departure from free kinetic transport.

The remainder of the paper is organized as follows. \cref{sec:kinetic_formulation} reviews the kinetic formulation and develops the characteristic and shock-supported representations underlying the proposed method. \cref{sec:ROM} presents the structure-informed \acp{ROM}, including kinetic lifting, characteristic pullback, shock-manifold registration, reduced shock and defect dynamics, and reconstruction. \cref{sec:numerical_examples} presents numerical examples for scalar conservation laws with shock formation and propagation. Conclusions are provided in \cref{sec:conclusion}.

\section{Kinetic formulation of scalar hyperbolic conservation laws}
\label{sec:kinetic_formulation}

We consider a scalar hyperbolic conservation law posed on $\Omega_{\bm x}\subseteq\mathbb{R}^{d}$,
\begin{equation}
    \partial_t u(\bm{x},t) + \nabla_{\bm{x}}\cdot \bm{f}\bigl(u(\bm{x},t)\bigr) = 0,
    \qquad
    (\bm{x},t)\in\Omega_{\bm x}\times(0,T],
    \label{eq:HCL}
\end{equation}
with initial condition
\begin{equation}
    u(\bm{x},0)=u_0(\bm{x}).
\end{equation}
Here, $u(\bm{x},t):\Omega_{\bm x}\times[0,T]\rightarrow\mathbb{R}$ denotes the conserved scalar quantity and $\bm{f}:\mathbb{R}\rightarrow\mathbb{R}^{d}$ is the flux function. For example, \(\bm f(u)=u^2/2\) corresponds to the inviscid Burgers equation, which is widely used as a canonical model of nonlinear wave propagation and shock formation, while \(\bm f(u)=\frac{u^2}{u^2+M(1-u)^2}\) is the Buckley--Leverett fractional-flow function, which is used to describe immiscible displacement in porous media, with \(M>0\) denoting the mobility ratio.

Even for smooth initial data and flux functions, nonlinear characteristics may intersect and form shocks in finite time. Since weak solutions are generally non-unique, an entropy condition is required to select the physically admissible solution~\cite{kruvzkov1970first,dafermos2005hyberbolic,leveque2002finite}.

\subsection{Kinetic entropy defect}
\label{subsec:KED}
The kinetic formulation lifts the scalar solution to an augmented space by introducing a kinetic coordinate $\xi\in\Omega_\xi\subset\mathbb{R}$. For each $(\bm{x},t)$, define the signed kinetic function
\begin{equation}
    \chi(\bm{x},t,\xi) := \chi\bigl(u(\bm{x},t),\xi\bigr) =
    \begin{cases}
        1,  & 0<\xi<u(\bm{x},t),\\[1mm]
       -1,  & u(\bm{x},t)<\xi<0,\\[1mm]
        0,  & \text{otherwise}.
    \end{cases}
    \label{eq:kinetic_function}
\end{equation}
For a nonnegative solution, $\chi$ reduces to the indicator of the kinetic interval $0<\xi<u(\bm{x},t)$. The physical state is recovered exactly through
\begin{equation}
    u(\bm{x},t) = \int_{\Omega_\xi}
    \chi(\bm{x},t,\xi)\,d\xi.
    \label{eq:kinetic_reconstruction}
\end{equation}
Equivalently, the unique entropy solution of~\cref{eq:HCL} is characterized by the transport equation
\begin{equation}
    \partial_t\chi + \bm{f}'(\xi)\cdot\nabla_{\bm{x}}\chi = \partial_\xi m,
    \label{eq:kinetic_equation}
\end{equation}
subject to the initial condition $\chi_0(\bm x, \xi)$. Here, the transport velocity $\bm{f}'(\xi)$ is determined entirely by the flux, and $m(\bm x, t, \xi)$ is referred to as the \ac{KED}, which has three structural properties (see~\cite{perthame2002kinetic} for details) that are central to the reduced-order formulation later:

\begin{property}[Nonnegativity and entropy admissibility]\label{prop:nonnegativity}
The nonnegativity of \(m\) encodes the entropy condition. Let  \(\phi:\mathbb{R}\rightarrow\mathbb{R}\) be a convex entropy and let \(\bm{f}_\phi\) denote its associated entropy flux, defined by
    \begin{equation}
    \bm{f}_\phi'(v)    =  \phi'(v)\bm{f}'(v).
    \end{equation}
    Testing~\cref{eq:kinetic_equation} with \(\phi'(\xi)\), integrating over \(\xi\), and assuming that the kinetic boundary terms vanish gives
    \begin{equation}
        \partial_t\phi(u)  + \nabla_{\bm{x}}\cdot\bm{f}_\phi(u) = -  \int_{\Omega_\xi}\phi''(\xi)\,m(\bm{x},t,d\xi) \leq 0.
    \end{equation}
    The inequality follows from the convexity of \(\phi\) and the nonnegativity of \(m\). Thus, a single kinetic measure represents the dissipation associated with every convex entropy pair, and \(m\) is constrained by the entropy condition rather than acting as an arbitrary closure term.
\end{property}

\begin{property}[Bounded support in the kinetic coordinate]\label{prop:bounded_support}
    If \(u_0\in L^\infty(\Omega_{\bm x})\), then
    \begin{equation}\label{eq:bounded_support}
        \operatorname{supp}_{\xi}m  \subseteq \Omega_\xi:= \left[ -\|u_0^-\|_{L^\infty(\Omega_{\bm x})}, \|u_0^+\|_{L^\infty(\Omega_{\bm x})}\right],
    \end{equation}
    where
    \begin{equation}
        u_0^+ = \max(u_0,0), \qquad u_0^- = \max(-u_0,0).
    \end{equation}
    Hence, only the kinetic levels attainable by the physical solution are relevant. This property allows the kinetic coordinate to be restricted to a finite interval in the numerical formulation.
\end{property}

\begin{property}[Localization in space--time]\label{prop:local}
    If \(u\) is continuously differentiable on an open region \(\mathcal{D}\subseteq\Omega_{\bm x}\times(0,T)\), then
    \begin{equation}
        m\big|_{\mathcal{D}\times\Omega_\xi}
        =0.
    \end{equation}
 The kinetic equation therefore reduces to homogeneous linear transport wherever the physical solution remains smooth. For piecewise smooth entropy solutions, the nontrivial defect is concentrated on the space-time shock set. At these locations, intersecting characteristics cause the breakdown of smooth transport, and \(\partial_\xi m\) supplies the entropy-producing correction required to recover the unique physically admissible solution.
\end{property}

The kinetic formulation thus lifts the nonlinear conservation law~\eqref{eq:HCL} into the higher-dimensional phase space \((\bm{x},t,\xi)\), where its dynamics~\eqref{eq:kinetic_equation} are separated into a known linear transport component and a localized, entropy-constrained defect component. The reduced-order formulation in~\cref{sec:ROM} exploits this structure by prescribing the kinetic transport velocity \(\bm{f}'(\xi)\) from the governing flux while learning the defect contribution from data.

\begin{remark}[Connection with level-set representations]
For each fixed kinetic coordinate \(\xi\), the function
\begin{equation}
    \ell(\bm{x},t,\xi) = u(\bm{x},t)-\xi
\end{equation}
defines the spatial level set \(\{\bm{x}: \ell(\bm{x},t,\xi)=0\}\). The kinetic function \(\chi(\bm{x},t,\xi)\) may therefore be interpreted as the collection of indicators of all level sets of \(u\), continuously parameterized by \(\xi\), rather than with a single interface. The use of level-set representations for \ac{ROM} of transport-dominated dynamics has also been explored, for example, in the appendix of Ref.~\cite{lu2020lagrangian}.
\end{remark}

\subsection{Characteristic representation}\label{subsec:charac}
For a fixed kinetic level $\xi$, the homogeneous part of \cref{eq:kinetic_equation} has characteristics
\begin{equation}
    \widehat{\bm{x}}(t;\bm{a},\xi) =    \bm{a}+\bm{f}'(\xi)t,
    \label{eq:charac}
\end{equation}
where $\bm{a}$ denotes the characteristic foot at $t=0$. Define the continuous pullback
\begin{equation}
    P(\bm{a},t,\xi) :=   \chi\left(  \bm{a}+\bm{f}'(\xi)t,t, \xi \right).
    \label{eq:con_pullback}
\end{equation}
Pulling back \cref{eq:kinetic_equation} along the free kinetic characteristics yields, in the sense of distributions,
\begin{equation}
    \partial_t P(\bm{a},t,\xi) =   \partial_\xi m \left( \bm{a}+\bm{f}'(\xi)t, t, \xi \right).
    \label{eq:pulled_back_kinetic}
\end{equation}
Thus, $P$ is constant along free kinetic characteristics wherever $m=0$, in particular in regions where the entropy solution is smooth. Any temporal change of the pulled-back kinetic state is caused by the \ac{KED}.

Integrating~\cref{eq:pulled_back_kinetic} from \(0\) to \(t\) along a characteristic and using \(P(\bm a,0,\xi)=\chi_0(\bm a,\xi)\) gives
\[
P(\bm a,t,\xi) = \chi_0(\bm a,\xi) + \int_0^t \partial_\xi m\bigl(\bm a+\bm f'(\xi)\tau,\tau,\xi\bigr)\,d\tau .
\]
Returning to physical coordinates with \(\bm a=\bm x-\bm f'(\xi)t\) yields the decomposition
\begin{equation}
\chi(\bm x,t,\xi)
=
\underbrace{
\chi_0\!\left(\bm x-\bm f'(\xi)t,\xi\right)
}_{\text{free transport}}
+
\underbrace{
\chi_{\mathrm{def}}(\bm x,t,\xi)
}_{\text{defect correction}}.
\end{equation}
where the first term is the kinetic state obtained by transporting the initial condition along the known free characteristics and the second correction term \(\chi_{\mathrm{def}} = \int_0^t \partial_\xi m\bigl(\bm x-\bm f'(\xi)(t-\tau),\tau,\xi\bigr)\,d\tau\) accumulates the kinetic defect encountered along the characteristic and therefore measures the departure of the entropy solution from free transport. This decomposition separates the known transport, determined entirely by the flux, from the defect-driven dynamics that remain to be modeled. For \ac{ROM}, the former need not be learned from data; the structure of the latter is examined next.

\subsection{Shock-supported structure of the KED}
\label{subsec:shock_supp}
The defect-driven contribution identified above is not generated throughout the full space--time domain. As discussed in~\cref{subsec:KED}, the KED vanishes in smooth regions and is localized to shock-producing events. This localization becomes more explicit when the entropy solution contains a sufficiently regular shock manifold. Following~\cite{srivastava2026computable}, let \(\Gamma\subset\Omega_{\bm x}\times(0, T)\) denote a sufficiently regular space-time shock set, whose spatial section at time $t$ is parameterized by
\begin{equation}
    \bm{s}(t,\bm{\theta}): (0,T)\times\Omega_{\bm \theta}\rightarrow\Omega_{\bm x},
    \qquad   \Omega_{\bm \theta}\subseteq\mathbb{R}^{d-1}.
    \label{eq:shock_par}
\end{equation}
Since the KED vanishes away from the shock set, a measure supported on this manifold can be represented formally as
\begin{equation}
    m(\bm{x},t,\xi) =    \int_{\Omega_{\bm \theta}}
    \overline m(t,\bm{\theta},\xi) \,    \delta\!\left(
        \bm{x}-\bm{s}(t,\bm{\theta})    \right) \,d\bm{\theta},
    \label{eq:shock_supported_ked}
\end{equation}
where \(\overline m(t,\bm{\theta},\xi)\) is the KED density restricted to the shock coordinates. Thus, although \(m\) is formally defined over the augmented space \((\bm{x},t,\xi)\), its nontrivial spatial dependence is confined to the lower-dimensional moving shock manifold.

Substitution of \cref{eq:shock_supported_ked} into the kinetic equation~\cref{eq:kinetic_equation} gives
\begin{equation}
    \partial_t\chi +    \bm{f}'(\xi)\cdot\nabla_{\bm{x}}\chi =    \int_{\Omega_{\bm \theta}}  \partial_\xi\overline m(t,\bm{\theta},\xi)    \,    \delta\!\left(  \bm{x}-\bm{s}(t,\bm{\theta})    \right)    \,d\bm{\theta}.
    \label{eq:shock_supported_kinetic_equation}
\end{equation}
Hence, the departure from free kinetic transport occurs only when characteristics interact with the shock manifold. In one spatial dimension, an isolated shock is described by a trajectory \(x=s(t)\), and \cref{eq:shock_supported_ked} reduces to

\begin{equation}
    m(x,t,\xi)  = \overline m(t,\xi)\,   \delta\bigl(x-s(t)\bigr).
\end{equation}
For multiple shocks, the measure may be represented by a corresponding sum of shock-supported contributions between interaction or topology-change events.

\cite{srivastava2026computable} exploits this shock-supported representation to constructively evaluate the KED from shock trajectories together with solution information. Here, the same localization is used instead as a model-reduction principle.  We do not require explicit evaluation of \(\overline m\); rather, the shock-supported structure identifies the dynamically active departure from free transport as a lower-dimensional moving event. This motivates extracting its geometry, registering the defect-driven evolution in shock-attached coordinates, and learning the resulting shock and defect dynamics directly from solution snapshots. In the data-driven construction below, we do not seek to recover $m$ itself; instead, temporal changes in the characteristic pullback provide an empirical approximation of the source $\partial_\xi m$ that drives departures from free transport.

\section{Data-driven ROMs via KED}\label{sec:ROM}
The kinetic formulation in \cref{sec:kinetic_formulation} separates the entropy solution into known free transport and a correction generated by the \ac{KED}. Moreover, \cref{subsec:shock_supp} shows that the source of this correction is concentrated on a lower-dimensional moving shock manifold. We construct a data-driven \ac{ROM} framework that uses both properties directly. Rather than reducing the full physical solution \(u(\bm{x},t)\) or the complete kinetic field, we prescribe the known characteristic transport analytically and approximate only the remaining defect-driven dynamics. Specifically, solution snapshots are lifted to kinetic space and pulled back along the free characteristics; temporal changes in the pullback are used to estimate the source \(\partial_\xi m\), and its moving support is registered in shock-attached coordinates. The resulting shock geometry and registered defect-driven source become our learning target whose dynamics are learned from data.

\subsection{Regularized kinetic lifting and decoding}\label{subsec:kinetic-lifting}
Let
\[
u^n(\bm{x}) := u(\bm{x},t_n), \qquad 0=t_0<t_1<\cdots<t_{N_t}=T,
\]
denote snapshots of the entropy solution of \cref{eq:HCL}. For the nonnegative solutions considered in the numerical examples,
$u^n(\bm{x})\geq 0$, and the kinetic function in~\cref{eq:kinetic_function} reduces to the subgraph indicator
\[
\chi^n(\bm{x},\xi) =\mathbf{1}_{\{0<\xi<u^n(\bm{x})\}}=H(\xi)H\!\left(u^n(\bm{x})-\xi\right).
\]
This product representation makes the two boundaries of the occupied kinetic interval explicit: $H(\xi)$ restricts the lift to positive kinetic levels, while $H(u^n-\xi)$ restricts it to levels below the physical state.

To obtain a numerically smooth kinetic field, we replace the discontinuous Heaviside function by \(H_\varepsilon(v) =\frac{1}{2}\left[1+\tanh\left(\frac{v}{\varepsilon}\right) \right]\), and define the diffuse kinetic lift used in the computations by
\begin{equation}
\psi_\varepsilon^n(\bm{x},\xi) = H_\varepsilon(\xi) H_\varepsilon\!\left(u^n(\bm{x})-\xi\right).
\label{eq:diffuse-kinetic-lift}
\end{equation}
For every fixed $(\bm{x},\xi)$ away from the boundaries $\xi=0$ and $\xi=u^n(\bm{x})$,
\[
\psi_\varepsilon^n(\bm{x},\xi) \longrightarrow \chi^n(\bm{x},\xi) \qquad\text{as }\varepsilon\rightarrow 0.
\]

Let \(\{\xi_j\}_{j=1}^{N_\xi}\subset \Omega_\xi\) denote the discrete kinetic grid, with associated quadrature weights $\{w_j\}_{j=1}^{N_\xi}$. By the maximum principle, $\Omega_\xi$ is bounded by 0 and \(\|u_0\|_{L^\infty(\Omega)}\) for nonnegative solutions. In the regularized computations, we use a slightly enlarged $\Omega_\xi$ to resolve the diffuse tails near the endpoints.

Based on~\cref{eq:kinetic_reconstruction}, the decoding of $u^n(\bm x)$ can be approximated by
\begin{equation}\label{eq:decode}
u^n(\bm x) \approx  \sum_{j=1}^{N_\xi}w_j\psi_\varepsilon^n(\bm{x},\xi_j).
\end{equation}
At finite $\varepsilon$ and finite kinetic resolution, the quadrature in~\cref{eq:decode} introduces a small deterministic bias. As a numerical remedy, we precompute its scalar input--output relation over the admissible range of $u$ and invert that relation by monotone interpolation when decoding the predicted kinetic state.

\begin{figure}[!h]
\centering
\includegraphics[height = 0.42\textwidth]{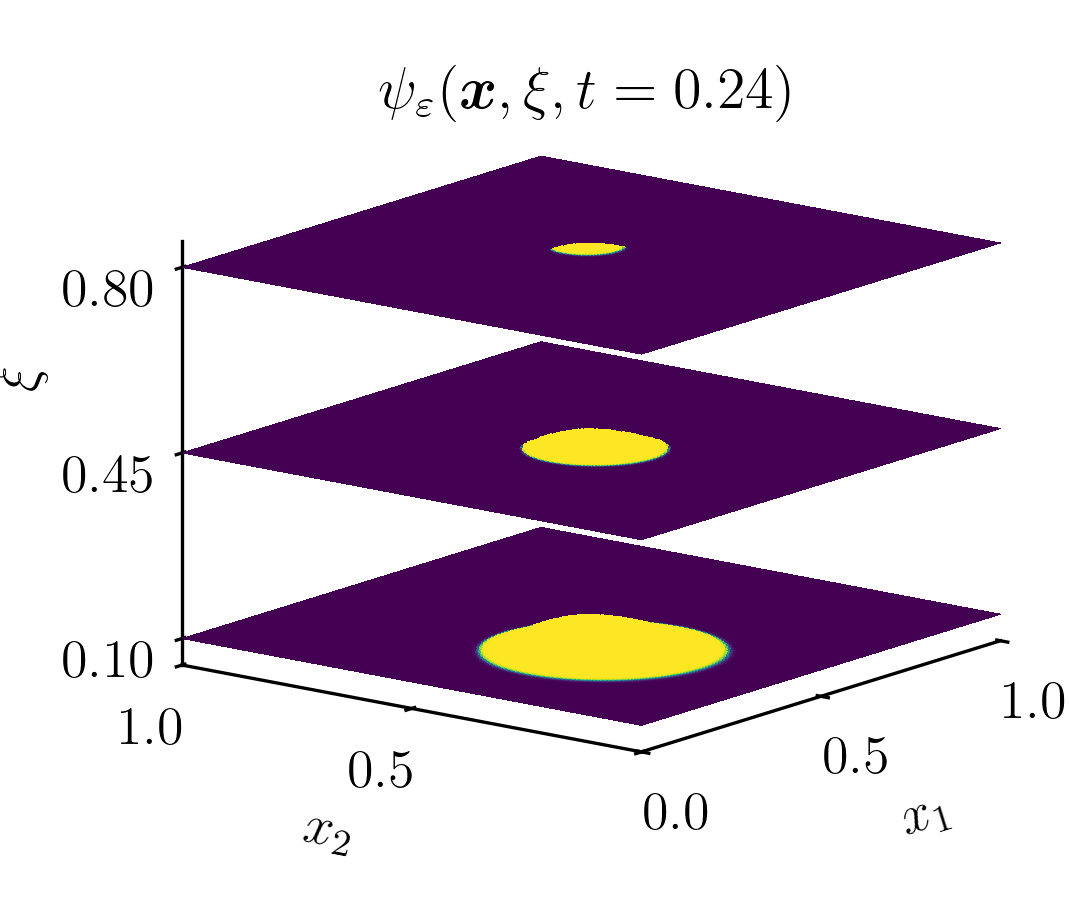}
\hfill
\includegraphics[height= 0.42\textwidth]{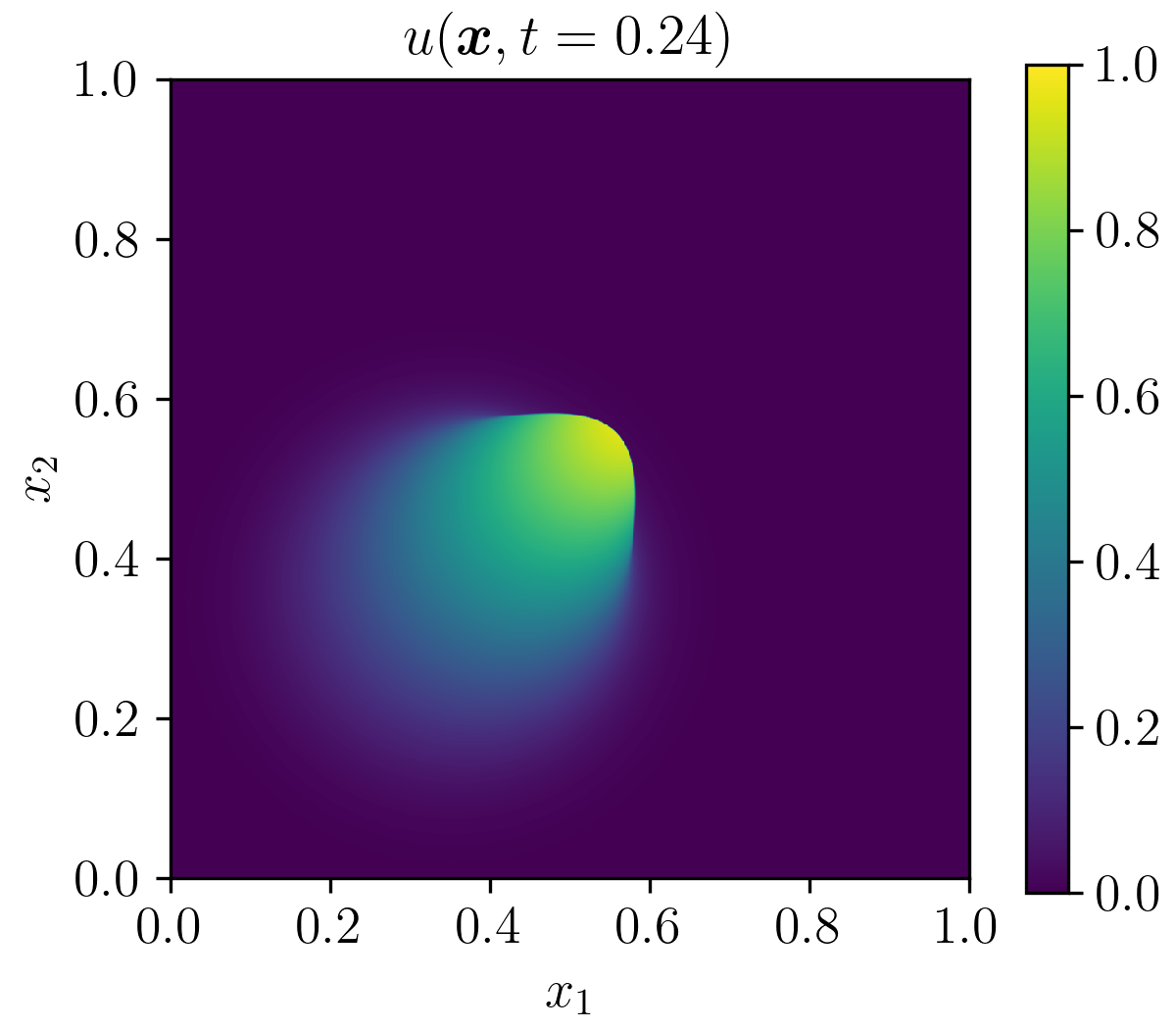}
\caption{Diffuse kinetic lifting and decoding at $t=0.24$ for numerical example in \cref{subsec:burgers_2d_gaussian}. Left: selected kinetic slices; Right: the decoded physical state. }
\label{fig:kinetic_lifting_defect}
\end{figure}

\cref{fig:kinetic_lifting_defect} illustrates the diffuse kinetic lift at selected values of $\xi$ and the physical state recovered by integration in the kinetic coordinate using~\cref{eq:decode}.

\subsection{Characteristic pullback and empirical defect-driven evolution}\label{subsec:cha-pullback}
For each kinetic level \(\xi_j\), the characteristic velocity \(\bm{f}'(\xi_j)\) is prescribed by the flux. Recall from \cref{subsec:charac} that the continuous characteristic pullback \cref{eq:con_pullback} satisfies \cref{eq:pulled_back_kinetic}. Thus, temporal variation along the known kinetic characteristics measures the departure from free transport generated by the \ac{KED}.

For the diffuse kinetic lift, we approximate this variation over the interval \([t^n,t^{n+1}]\) and associate it with the temporal midpoint \(t^{n+1/2}=(t^n+t^{n+1})/2\). Rather than explicitly constructing the complete pullback on a characteristic-coordinate grid, we evaluate the characteristic difference directly at a desired midpoint location \(\bm x_k\). Specifically, we define 
\begin{equation}\label{eq:char-diff}
G_\varepsilon^{n+1/2}(\bm x_k,\xi_j)
=
\frac{
\psi_\varepsilon^{n+1}
\left(
\bm x_k+\frac{\Delta t}{2}\bm f'(\xi_j),\xi_j
\right)
-
\psi_\varepsilon^{n}
\left(
\bm x_k-\frac{\Delta t}{2}\bm f'(\xi_j),\xi_j
\right)
}{\Delta t}.
\end{equation}
The off-grid values of \(\psi_\varepsilon^n\) are obtained by spatial interpolation. This expression is the centered finite-difference approximation of the temporal derivative of the characteristic pullback and therefore serves as an empirical regularized approximation of \(\partial_\xi m\) at \((\mathbf x_k,t^{n+1/2},\xi_j)\).


Before shock formation, the exact characteristic pullback is stationary. For the diffuse and discretized kinetic field,  \(G_\varepsilon^{n+1/2}\) therefore remains small. Let \(n_{\mathrm{ev}}\) denote the first event index, selected such that $t^{n_{\mathrm{ev}}-1}$ is near the observed or known breaking time. The shock registration and reduced shock and defect models are constructed from the post-breaking midpoint snapshots associated with,
\begin{equation}
  t^{n+1/2},\quad n\in \mathcal{N}_{\mathrm{ev}} = \{n_{\mathrm{ev}}, \ldots, N_{\mathrm{tr}}\}.
    \label{eq:event_training_indices}
\end{equation}
The pre-event evolution is retained during reconstruction but is governed by homogeneous characteristic transport rather than being included as nearly zero data in the defect reduction.

\subsection{Shock-manifold extraction and registration}\label{subsec:event_registration}
Although the characteristic difference in \cref{eq:char-diff} removes the free transport generated by \(\bm f'(\xi)\), the support of \(G_\varepsilon^{n+1/2}\) remains attached to an evolving shock manifold. Direct learning of this field would therefore mix shock translation, front deformation, and changes in the localized defect shape, making it extremely challenging for ROM constructions.  Motivated by the structural study in~\cref{subsec:shock_supp}, we use registration to represent these contributions in a common shock-attached coordinate system. The construction produces two quantities from \(G_\varepsilon^{n+1/2}(\bm x, \xi)\): (i) a physical embedding \(\bm c^{n+1/2}(\bm\theta)\) of the shock manifold and (ii) the registered defect-driven source \(\tilde G_\varepsilon^{n+1/2}(\bm\theta,\eta,\xi)\). Here \(\bm\theta\in\Omega_{\bm \theta}\subset\mathbb R^{d-1}\) is a tangential coordinate on a fixed reference domain, and \(\eta\in\Omega_\eta\subset\mathbb R\) is a signed shock-normal coordinate. For simplicity, one may take \(\Omega_{\bm \theta}=[0,1]^{d-1}\), with the appropriate boundary identifications for closed fronts. 

\subsubsection{Defect-based shock-manifold embedding}\label{subsubsec:shock}
The shock embedding provides a data-derived approximation of the shock parameterization introduced in \cref{eq:shock_par}:
\begin{equation}
    \bm c^{n+1/2}(\bm\theta)\approx\bm s(t^{n+1/2},\bm\theta),\qquad\bm\theta\in\Omega_{\bm \theta}=[0,1]^{d-1}.
    \label{eq:data_shock_embedding}
\end{equation}
It is constructed by selecting resolved kinetic levels, combining their defect information in physical coordinates, and extracting the resulting shock ridge.

For every kinetic level \(\xi_j\), define the defect energy
\begin{equation}
    e_j^{n+1/2}= \left(\int_{\Omega_{\bm x}} \left|G_\varepsilon^{n+1/2}(\bm x,\xi_j) \right|^2 \,\mathrm d\bm x\right)^{1/2}.
    \label{eq:defect_level_energy}
\end{equation}
The active kinetic set is
\begin{equation}
    \mathcal J^{n+1/2}=\left\{ j:e_j^{n+1/2}> \kappa_{\mathrm{act}}\max_{\widetilde j} e_{\widetilde j}^{n+1/2}\right\}, \quad 0< \kappa_{\mathrm{act}}<1.
    \label{eq:active_kinetic_set}
\end{equation}
Only these levels are used to infer the shock geometry. Because the same physical shock generates defect activity across multiple kinetic levels, we aggregate the active slices to obtain a scalar indicator of the shock location, i.e., the aggregate event density
\begin{equation}
    R^{n+1/2}(\bm x) = \sum_{j\in\mathcal J^{n+1/2}}w_j \left| G_\varepsilon^{n+1/2} \left( \bm x, \xi_j\right)\right|,
    \label{eq:physical_event_density}
\end{equation}
where \(w_j\) are the kinetic quadrature weights.

Within a local graph chart, choose an origin \(\bm x_o\), a tangential basis \(\bm T\in\mathbb R^{d\times(d-1)}\), and a transverse unit vector \(\bm d\in\mathbb R^d\) satisfying
\begin{equation}
    \bm T^T\bm T=\bm I, \qquad  \bm T^T\bm d=\bm 0, \qquad \|\bm d\|_2=1.
    \label{eq:local_shock_frame}
\end{equation}
The associated coordinates are
\begin{equation}
    \bm x(\bm q,p) =\bm x_o+\bm T\bm q+p\bm d, \qquad (\bm q,p)\in\Omega_{\bm q}^{\mathrm{probe}}\times\Omega_{p}\subset \mathbb R^{d-1}\times \mathbb R.
    \label{eq:local_shock_coordinates}
\end{equation}
Here, \(\Omega_{\bm q}^{\mathrm{probe}}\) is the prescribed tangential search domain, while \(\Omega_{p}\) is the transverse interval over which the defect is sampled. The resolved tangential support is subsequently identified as a subset of \(\Omega_{\bm q}^{\mathrm{probe}}\):
\begin{equation}
    \mathcal Q^{n+1/2} = \left\{ \bm q\in\Omega_{\bm q}^{\mathrm{probe}}:E_\Gamma^{n+1/2}(\bm q) > \kappa_\Gamma \max_{\widetilde{\bm q}} E_\Gamma^{n+1/2}(\widetilde{\bm q}) \right\}, 0<\kappa_\Gamma<1,
    \label{eq:resolved_tangential_support}
\end{equation}
where the tangential defect energy
\[
    E_\Gamma^{n+1/2}(\bm q) = \left( \sum_{j\in\mathcal J^{n+1/2}}w_j\int_{\Omega_{p}}\left|G_\varepsilon^{n+1/2}\left( \bm x(\bm q,p), \xi_j \right)\right|^2 \,\mathrm dp \right)^{1/2},
\]
is used to determine the tangential extent over which the ridge is resolved.

For every $\bm q\in\mathcal Q^{n+1/2}$, the raw transverse shock position is defined by the strongest defect response,
\begin{equation}
    p_{\mathrm{raw}}^{n+1/2}(\bm q) = \underset{p\in\Omega_{p}}{\operatorname{arg\,max}}\; R^{n+1/2}\bigl(\bm x(\bm q,p)\bigr).
    \label{eq:raw_shock_manifold}
\end{equation}
On a discrete probe grid, the maximizer may be refined by local interpolation to obtain a subgrid approximation of $p_{\mathrm{raw}}^{n+1/2}(\bm q)$. A spatial fit or interpolation of the resulting raw ridge is then denoted by $p_\Gamma^{n+1/2}(\bm q)$. The particular refinement and spatial fitting procedures are problem-dependent and are specified in the numerical examples.

To represent shock manifolds with different tangential extents on the fixed reference domain \(\Omega_{\bm \theta}=[0,1]^{d-1}\), define the coordinatewise limits of the resolved
tangential support by
\begin{equation}
    q_r^-  =  \min_{\bm q\in\mathcal Q^{n+1/2}}q_r, \qquad   q_r^+ = \max_{\bm q\in\mathcal Q^{n+1/2}}q_r, \qquad r=1,\ldots,d-1.
    \label{eq:tangential_support_limits}
\end{equation}
These limits determine the coordinate-aligned chart domain
\begin{equation}
    \widehat{\mathcal Q}^{n+1/2} = \prod_{r=1}^{d-1}[q_r^-,q_r^+],
    \label{eq:tangential_chart_domain}
\end{equation}
which contains the resolved support \(\mathcal Q^{n+1/2}\). We then introduce the tangential registration map
\begin{equation}
    \bm\Phi^{n+1/2}: \Omega_{\bm \theta} \longrightarrow\widehat{\mathcal Q}^{n+1/2},
    \label{eq:tangential_registration_map}
\end{equation}
defined component-wise by
\begin{equation}
    \Phi_r^{n+1/2}(\bm\theta)  = q_r^-  + \theta_r\left(q_r^+-q_r^-\right),  \qquad  r=1,\ldots,d-1.
    \label{eq:affine_tangential_registration}
\end{equation}
Thus, \(\bm\theta\) specifies a normalized tangential location on a fixed reference domain, while \(\bm\Phi^{n+1/2}(\bm\theta)\) accounts for the time-dependent physical
extent of the shock chart. In general, \(\widehat{\mathcal Q}^{n+1/2}\) may contain points outside the resolved support \(\mathcal Q^{n+1/2}\); these points are excluded by the resolved support indicator when evaluating the registered defect. When \(d=2\), the tangential support is an interval and \(\widehat{\mathcal Q}^{n+1/2}=\mathcal Q^{n+1/2}\).

\begin{figure}[!h]
\centering
\includegraphics[width= 0.49\textwidth,trim={0 10 0 10},clip]{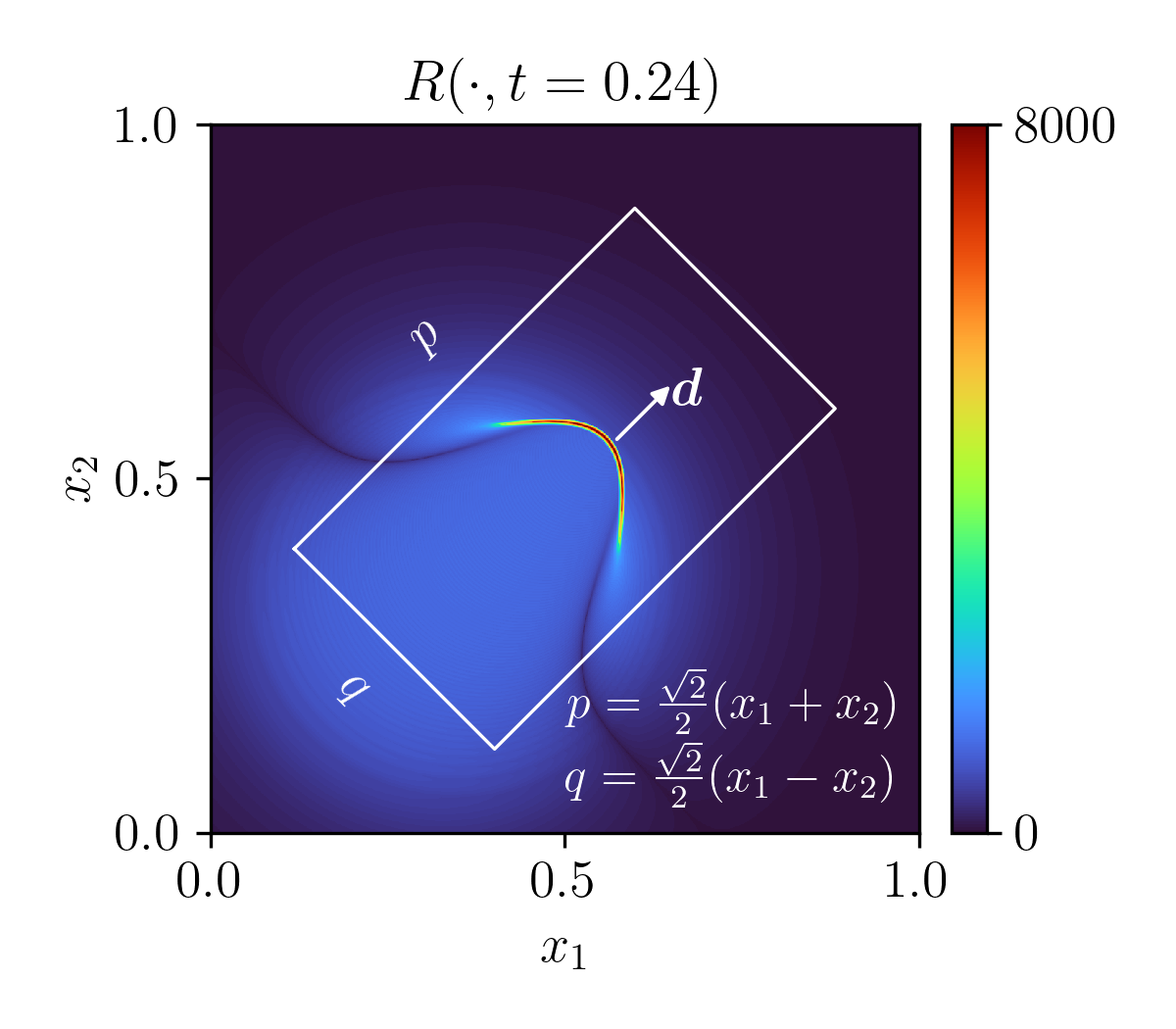}
\hfill
\includegraphics[width= 0.49\textwidth,trim={0 10 0 10},clip]{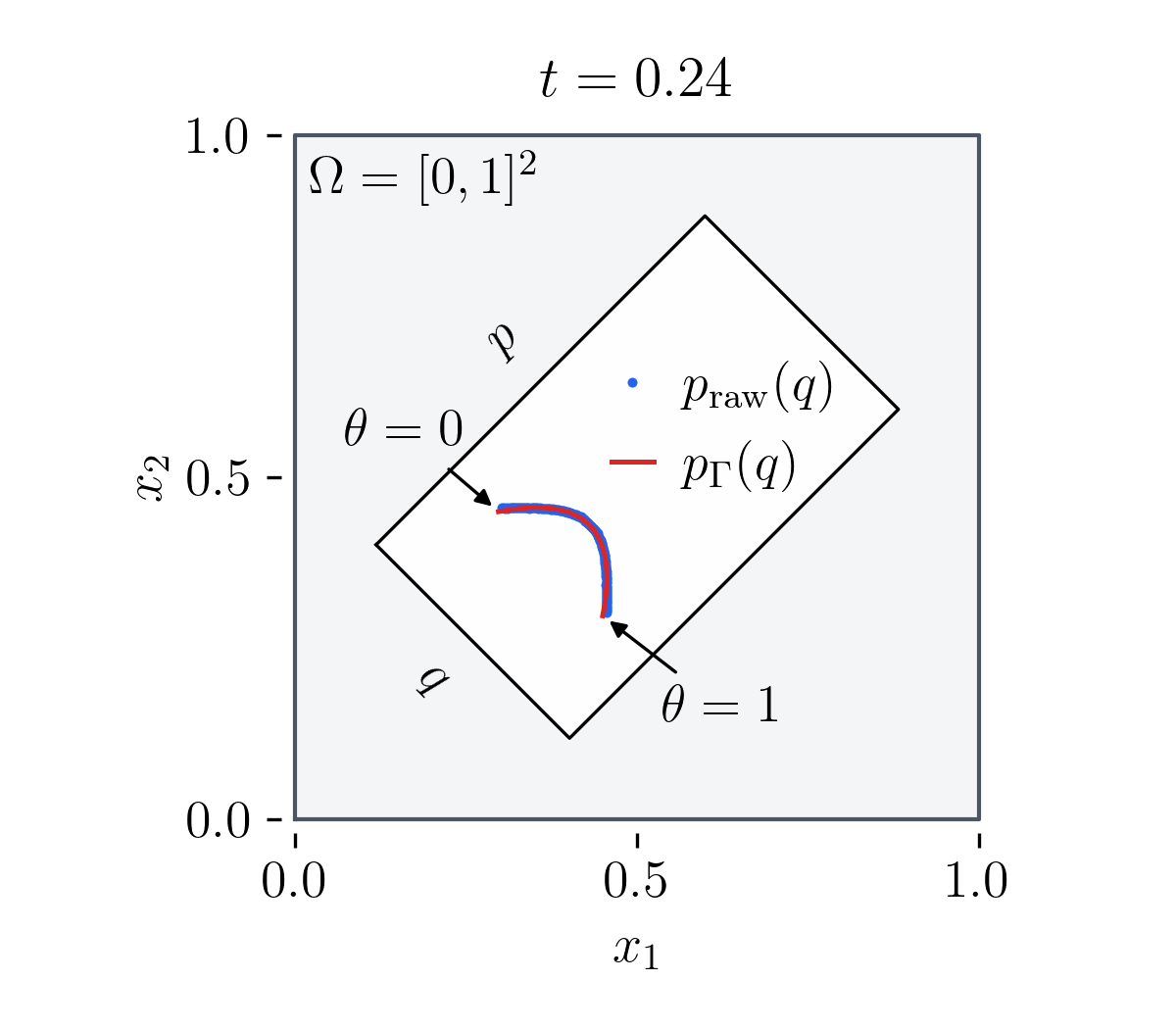}
\caption{Shock-manifold extraction and registration at $t=0.24$ for the two-dimensional example in \cref{subsec:burgers_2d_gaussian}. Left:  the aggregate event density; Right: the raw ridge $p_{\mathrm{raw}}(q)$, its fitted representation $p_\Gamma(q)$, and the corresponding embedding in the physical domain.}
\label{fig:shock}
\end{figure}

The shock embedding is then
\begin{equation}
    \bm c^{n+1/2}(\bm\theta)  = \bm x\left( \bm\Phi^{n+1/2}(\bm\theta), p_\Gamma^{n+1/2}  \left(\bm\Phi^{n+1/2}(\bm\theta)\right) \right),  \qquad \bm\theta\in\Omega_{\bm \theta}.
    \label{eq:shock_manifold_embedding}
\end{equation}
The same reference coordinate \(\bm\theta\) is subsequently used to represent the registered defect-driven source. Hence, \(\bm c^{n+1/2}\) records both the physical location and the tangential extent of the shock manifold on a fixed reference domain. \cref{fig:shock} illustrates the construction for the two-dimensional example.

\subsubsection{Registered defect-driven source}\label{subsubsec:registered_defect}
The registered defect-driven source is obtained by expressing the empirical defect-driven evolution in shock-attached coordinates. Define
\begin{equation}
    \bm X^{n+1/2}(\bm\theta,\eta) = \bm c^{n+1/2}(\bm\theta) + \eta\bm d,\qquad(\bm\theta,\eta)\in\Omega_{\bm \theta}\times\Omega_\eta,
    \label{eq:shock_transverse_chart}
\end{equation}
where \(\bm d\) is the transverse chart direction from \cref{eq:local_shock_frame}. 

\begin{figure}[!h]
\centering
\includegraphics[width= \textwidth,trim={0 20 0 20},clip]{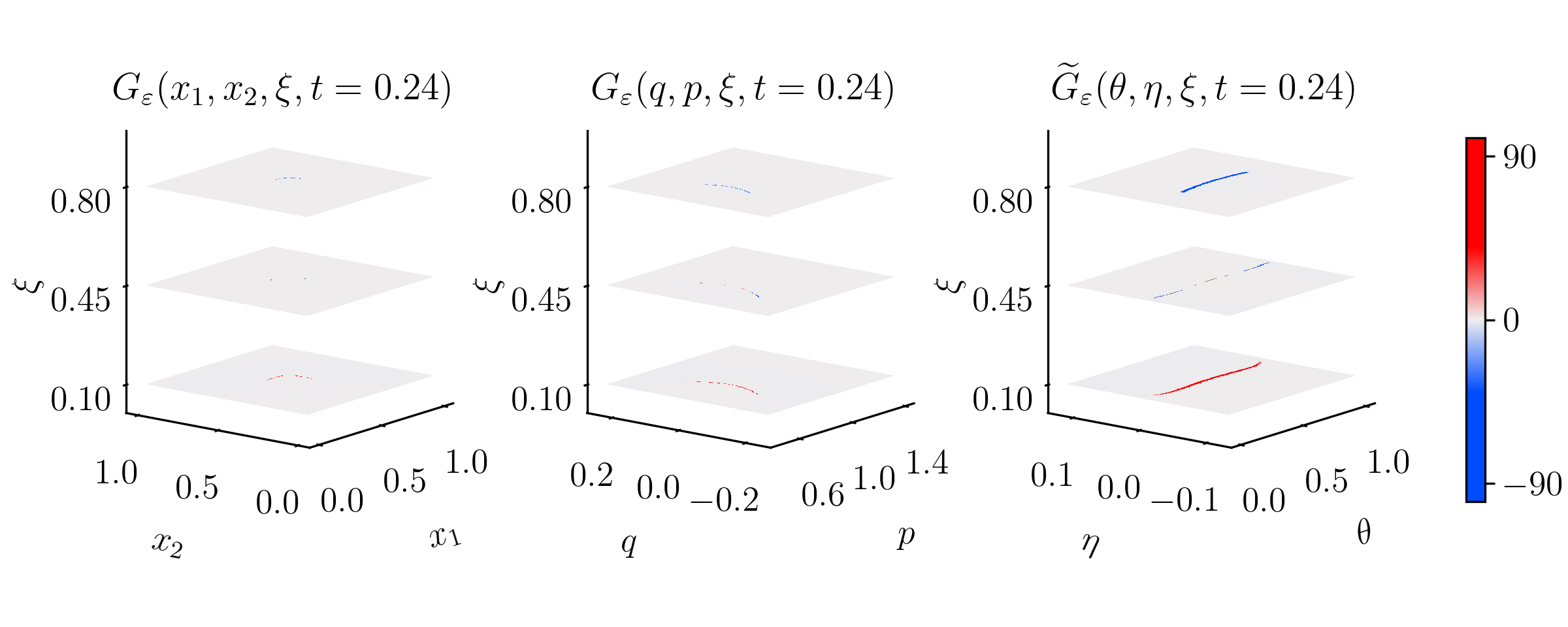}
\caption{Selected defect-source slices at $t=0.24$ for the two-dimensional example in \cref{subsec:burgers_2d_gaussian}. Left: physical coordinates \((x_1,x_2)\); Middle: graph-chart coordinates \((p,q)\); Right: registered coordinates \((\theta, \eta)\).}
\label{fig:defect}
\end{figure}

The registered defect-driven source is
\begin{equation}
    \widetilde G_\varepsilon^{n+1/2} (\bm\theta,\eta,\xi) = G_\varepsilon^{n+1/2} \left( \bm X^{n+1/2}(\bm\theta,\eta), \xi\right).
    \label{eq:registered_empirical_defect}
\end{equation}
Because \(G_\varepsilon^{n+1/2}\) is evaluated directly at physical midpoint locations, the registered source can be constructed without forming the complete defect field on the physical domain. For each registered point \(\bm X^{n+1/2}(\bm\theta,\eta)\), the characteristic difference in \cref{eq:char-diff} is evaluated directly from the two neighboring kinetic snapshots. This pointwise construction is particularly useful in multiple spatial dimensions, where assembling \(G_\varepsilon^{n+1/2}(\bm x,\xi)\) over the full physical-kinetic domain would be unnecessarily expensive. \cref{fig:defect} shows an example of the defect-driven source in different coordinate representations.

\subsubsection{\acp{ROM} for the shock geometry and registered defect-driven source}
\label{subsubsec:event_rom}

The shock embedding and registered defect-driven source are modeled by two
independent \acp{ROM} via \ac{DMD}. On the grid
$\{\bm\theta_k\}_{k=1}^{N_\theta}$, define the discrete shock-geometry state
\begin{equation}
    \bm z_c^{n+1/2}
    =
    \operatorname{vec}_{\bm\theta}
    \left[
        \bm c^{n+1/2}(\bm\theta_k)
    \right]
    \in \mathbb R^{dN_\theta}.
    \label{eq:discrete_shock_state}
\end{equation}
On the grids $\{\bm\theta_k\}_{k=1}^{N_\theta}$,
$\{\eta_\ell\}_{\ell=1}^{N_\eta}$, and
$\{\xi_j\}_{j=1}^{N_\xi}$, define the discrete registered defect state
\begin{equation}
    \bm z_G^{n+1/2}
    =
    \operatorname{vec}_{\xi,\bm\theta,\eta}
    \left[
        \widetilde G_\varepsilon^{n+1/2}
        (\bm\theta_k,\eta_\ell,\xi_j)
    \right]
    \in
    \mathbb R^{N_\xi N_\theta N_\eta}.
    \label{eq:discrete_registered_event_state}
\end{equation}

The same reduction and identification procedure is applied independently
to $\bm z_c$ and $\bm z_G$. For either state $\bm z_s$, $s\in\{c,G\}$,
the training snapshots are centered by their temporal mean and scaled by
a common scalar normalization. Let $\widehat{\bm z}_s^{n+1/2}$ denote
the resulting normalized snapshots and form
\begin{equation}
    \bm Z_{s,1}
    =
    \begin{bmatrix}
        \widehat{\bm z}_s^{\,n_{\mathrm{ev}}+1/2}
        & \cdots &
        \widehat{\bm z}_s^{\,N_{\mathrm{tr}}-1/2}
    \end{bmatrix},
    \quad
    \bm Z_{s,2}
    =
    \begin{bmatrix}
        \widehat{\bm z}_s^{\,n_{\mathrm{ev}}+3/2}
        & \cdots &
        \widehat{\bm z}_s^{\,N_{\mathrm{tr}}+1/2}
    \end{bmatrix}.
    \label{eq:dmd_snapshot_matrices}
\end{equation}
Using the rank-$r_s$ truncated singular value decomposition
$\bm Z_{s,1}\approx
\bm U_{s,r_s}\bm\Sigma_{s,r_s}\bm V_{s,r_s}^{\top}$,
the corresponding reduced coordinates are
\begin{equation}
    \bm z_{s,r}^{n+1/2}
    =
    \bm U_{s,r_s}^{\top}
    \widehat{\bm z}_s^{n+1/2}.
\end{equation}
Each reduced state is then evolved by its own affine \ac{DMD} model,
\begin{equation}
    \bm z_{s,r}^{n+3/2}
    =
    \bm A_{s,r}\bm z_{s,r}^{n+1/2}
    +
    \bm b_{s,r},
    \qquad
    s\in\{c,G\},\qquad n = 0, \cdots, N_\text{tr},\cdots,
    \label{eq:affine_dmd}
\end{equation}
where $(\bm A_{c,r},\bm b_{c,r})$ and
$(\bm A_{G,r},\bm b_{G,r})$ are identified independently from their
respective consecutive reduced snapshots. The ranks $r_c$ and $r_G$
are likewise selected independently from the corresponding singular-value
decays; the criterion used in the numerical examples is specified in
\cref{sec:numerical_examples}.

\begin{remark}\label{rmk:1d}
In one spatial dimension, the shock manifold at a fixed time consists
of isolated points. For a single shock, $\Omega_\theta$ is therefore
absent, and the shock-geometry state reduces to the scalar shock location
$c^{n+1/2}$. Rather than constructing a \ac{DMD} model for this scalar
quantity, its temporal evolution is approximated directly by a polynomial
fit over the training interval. The registered defect source depends only
on $(\eta,\xi)$, with discrete state
\begin{equation}
    \bm z_G^{n+1/2}
    =
    \operatorname{vec}_{\xi,\eta}
    \left[
        \widetilde G_\varepsilon^{n+1/2}(\eta_\ell,\xi_j)
    \right].
\end{equation}
Its reduced dynamics are evolved by the affine \ac{DMD} model described
above. Thus, in one spatial dimension the shock location is modeled
directly by polynomial regression, while the registered defect retains
the same reduced \ac{DMD} construction used in multiple spatial
dimensions.
\end{remark}

\subsubsection{Inverse registration}
\label{subsubsec:discrete_event_state}
During online prediction, the shock-geometry model provides the predicted embedding \(\bm c_{\mathrm{ROM}}^{n+1/2}(\bm\theta)\), while the registered-defect \ac{DMD} model provides \(\widetilde G_{\varepsilon,\mathrm{ROM}}^{n+1/2}\). The latter is defined on the fixed registered grid \(\Omega_{\bm\theta}\times\Omega_\eta\times\Omega_\xi\) and must therefore be inverse-registered to physical midpoint coordinates before it can be used in the kinetic update.

\begin{figure}[!h]
\centering
\includegraphics[width= \textwidth,trim={0 30 0 50},clip]{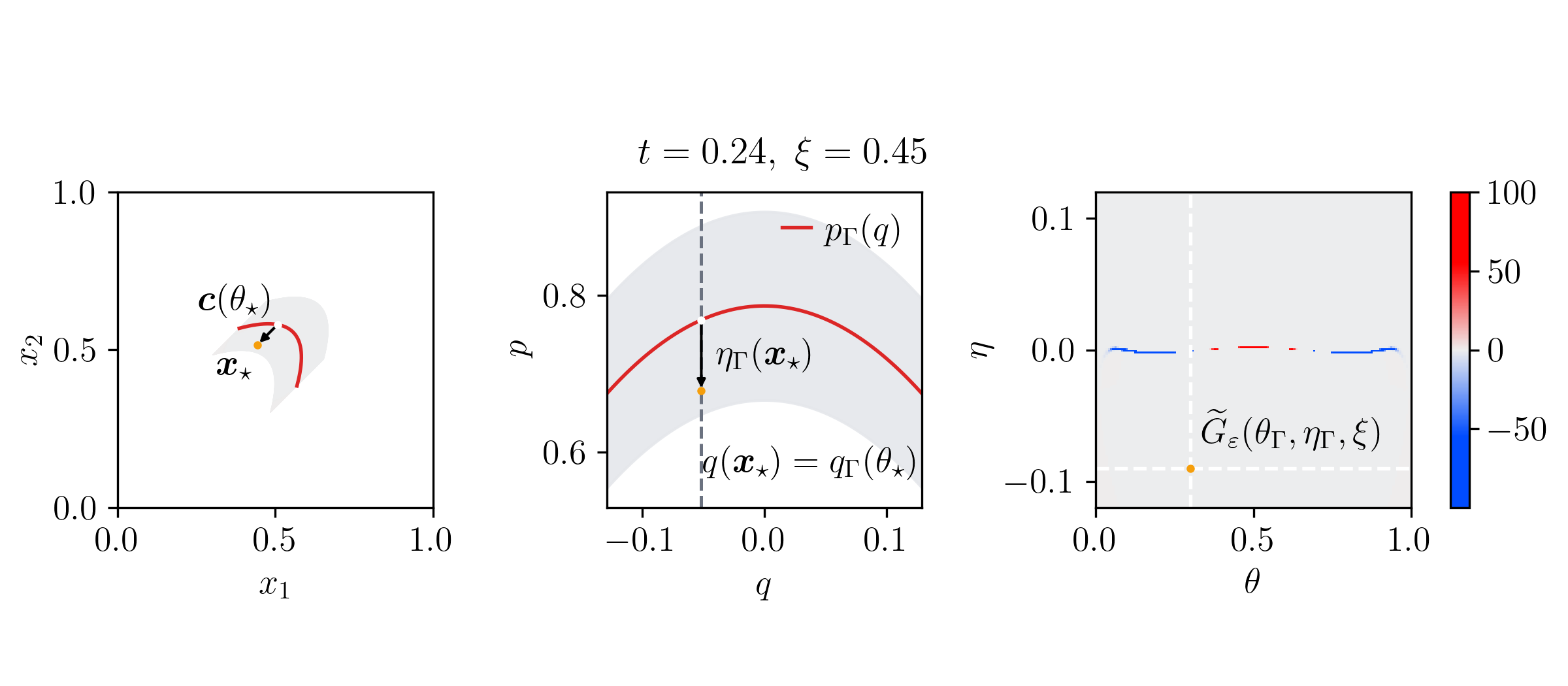}
\caption{Inverse registration at $t=0.24$ and $\xi=0.45$. From left to right, the physical point $\bm x_\star$ is mapped to the shock-fitted coordinates $(\theta_\star,\eta_\Gamma(\bm x_\star))$, where the registered defect $\widetilde G_\varepsilon$ is evaluated.}
\label{fig:inverse_registration_workflow}
\end{figure}

For every kinetic level \(\xi_j\), the inverse-registration workflow (see \cref{fig:inverse_registration_workflow}) first expresses the physical point \(\bm x_\star\) in the fixed graph-chart coordinates:
\begin{equation}
    \bm q(\bm x_\star) =   \bm T^T(\bm x_\star-\bm x_o), \qquad p(\bm x_\star)  =\bm d^T(\bm x_\star-\bm x_o).
    \label{eq:inverse_graph_coordinates}
\end{equation}
The ROM-predicted shock embedding is expressed in the same coordinates:
\begin{equation}
    \bm q_\Gamma^{n+1/2}(\bm\theta) = \bm T^T \left(  \bm c_{\mathrm{ROM}}^{n+1/2}(\bm\theta)-\bm x_o \right),\ p_\Gamma^{n+1/2}(\bm\theta) = \bm d^T \left( \bm c_{\mathrm{ROM}}^{n+1/2}(\bm\theta)-\bm x_o \right).
    \label{eq:reconstructed_graph_coordinates}
\end{equation}

Assuming that \(\bm q_\Gamma^{n+1/2}:\Omega_{\bm \theta}\rightarrow\mathbb R^{d-1}\) is one-to-one within the graph chart, we numerically recover its inverse through
\begin{equation}
    \bm\theta_\Gamma^{n+1/2}(\bm x_\star) = \underset{\tilde{\bm\theta}\in\Omega_{\bm \theta}}{\operatorname*{argmin}} \left\| \bm q(\bm x_\star) - \bm q_\Gamma^{n+1/2}(\tilde{\bm\theta})  \right\|_2^2.
    \label{eq:inverse_tangential_coordinate}
\end{equation}
The corresponding signed transverse displacement from the reconstructed shock is
\begin{equation}
    \eta_\Gamma^{n+1/2}(\bm x_\star) = p(\bm x_\star) -  p_\Gamma^{n+1/2} \left(  \bm\theta_\Gamma^{n+1/2}(\bm x_\star) \right).
    \label{eq:inverse_transverse_coordinate}
\end{equation}
Hence, \(\bm\theta_\Gamma^{n+1/2}(\bm x_\star)\) identifies the location along the shock chart, while \(\eta_\Gamma^{n+1/2}(\bm x_\star)\) gives the displacement from the reconstructed shock along the fixed transverse direction \(\bm d\).

The physical midpoint defect-driven source is then reconstructed by interpolation on the registered grid:
\begin{equation}
G_{\varepsilon,\mathrm{ROM}}^{n+1/2}(\bm x_\star,\xi)  =\widetilde G_{\varepsilon,\mathrm{ROM}}^{n+1/2} \left(\bm\theta_\Gamma^{n+1/2}(\bm x_\star),   \eta_\Gamma^{n+1/2}(\bm x_\star),   \xi \right) \mathbf 1_{\mathcal C_{\mathrm{ROM}}^{n+1/2}}(\bm x_\star).
    \label{eq:inverse_defect_registration}
\end{equation}
Here, \(\mathcal C_{\mathrm{ROM}}^{n+1/2}\) is the physical region covered by the reconstructed chart, characterized by

\begin{equation}
    \bm q(\bm x_\star) \in  \bm q_\Gamma^{n+1/2}(\Omega_{\bm \theta}),  \qquad  \eta_\Gamma^{n+1/2}(\bm x_\star)\in\Omega_\eta.
    \label{eq:inverse_registration_validity}
\end{equation}
The reconstructed defect is set to zero outside this region.

Finally, the physical kinetic update can be written as
\begin{equation}    \label{eq:physical_kinetic_update}
    \psi_{\varepsilon,\mathrm{ROM}}^{n+1}(\bm x_\star,\xi) =  \psi_{\varepsilon,\mathrm{ROM}}^n \left(  \bm x_\star-\Delta t\,\bm f'(\xi),     \xi  \right) 
    + \Delta t\,G_{\varepsilon,\mathrm{ROM}}^{n+1/2} \left( \bm x_\star-\frac{\Delta t}{2}\bm f'(\xi),  \xi \right).
\end{equation}
\cref{eq:physical_kinetic_update} is a characteristic semi-Lagrangian update of the diffuse kinetic state: the homogeneous contribution is obtained by tracing \(\bm x_\star\) backward along the analytically known kinetic characteristic over one time step, while the ROM supplies the defect-driven source at the corresponding characteristic midpoint. 

\Cref{fig:inverse_registration_workflow} illustrates this pointwise procedure for the two-dimensional example. A physical midpoint is first expressed in the graph coordinates \((\bm q,p)\), its registered coordinates \((\bm\theta_\Gamma,\eta_\Gamma)\) are recovered relative to the predicted shock, and the registered defect is then evaluated at those coordinates. The procedure is a point-wise coordinate lookup rather than a global deformation of the registered array. 

\subsection{Computational workflow}
\label{subsec:alg}

The complete offline-online procedure is summarized in \cref{alg:kd_rom}.
During the offline stage, solution snapshots are lifted to kinetic space, and
the empirical defect-driven source is evaluated along the known
characteristics. The shock manifold is then extracted, the defect-driven
source is registered in shock-attached coordinates, and separate reduced
models are constructed for the shock geometry and registered defect-driven
source. During online prediction, these models are evolved independently,
and the predicted shock geometry is used to inverse-register the predicted
defect-driven source. The resulting physical source is combined with the
prescribed characteristic transport to advance the kinetic state. The
graph-chart specifications and discretization parameters are provided with
the numerical examples in \cref{sec:numerical_examples}.

\begin{algorithm}[!htbp]
\caption{Kinetic-defect reduced-order model}
\label{alg:kd_rom}
\small

\textbf{Input:}
snapshots $\{u_\text{ref}^n\}_{n=1}^{N_\text{tr}}$ and the physical,
kinetic, and registration grids.

\medskip
\textbf{Offline stage:}
\begin{enumerate}
    \setlength{\itemsep}{0.15em}

    \item Construct the diffuse kinetic lift and empirical defect-driven
    source using
    \cref{eq:diffuse-kinetic-lift,eq:char-diff}.

    \item Extract the shock embedding $\bm c^{n+1/2}$ using
    \crefrange{eq:defect_level_energy}{eq:shock_manifold_embedding},
    and construct the registered defect-driven source using
    \cref{eq:shock_transverse_chart,eq:registered_empirical_defect}.

    \item Construct independent \acp{ROM} \eqref{eq:affine_dmd} for the shock geometry and
    registered defect-driven source as described in
    \cref{subsubsec:event_rom}.
\end{enumerate}

\medskip
\textbf{Online stage:}
\begin{enumerate}
    \setlength{\itemsep}{0.15em}

    \item Evolve the shock-geometry and registered-defect models \eqref{eq:affine_dmd} to obtain
    $\bm c_\text{ROM}^{n+1/2}$ and
    $\widetilde G_{\varepsilon,\text{ROM}}^{n+1/2}$.

    \item Inverse-register the predicted defect-driven source using
    \crefrange{eq:inverse_graph_coordinates}{eq:inverse_registration_validity}.

    \item Advance the kinetic field using
    \cref{eq:physical_kinetic_update} and decode the physical solution
    using \cref{eq:decode}.
\end{enumerate}

\textbf{Output:}
reconstructions $\{u_\text{ROM}^n\}_{n=1}^{N_\text{tr}}$ and predictions
$\{u_\text{ROM}^n\}_{n=N_\text{tr}+1}^{\ldots}$.
\end{algorithm}

\section{Numerical examples}
\label{sec:numerical_examples}
We assess the proposed kinetic-defect ROM on four scalar hyperbolic conservation-law problems of increasing complexity. The examples include the Burgers and Buckley-Leverett fluxes in one and two spatial dimensions and are chosen to examine shock formation, propagation, and multidimensional front evolution. In each case, the ROM is trained on a finite post-breaking time interval and subsequently evolved beyond the training window to assess its predictive capability.


For all examples, the ROMs are constructed following
\cref{alg:kd_rom}, and each DMD rank is selected independently as the
smallest rank retaining \(99\%\) of the corresponding snapshot energy.
We fix the diffuse-interface parameter at \(\varepsilon=0.01\). Since
\(H_\varepsilon\) varies on the kinetic scale \(\varepsilon\), the kinetic
grids are chosen with \(\Delta\xi=\varepsilon\), and \(\Omega_\xi\) is
extended beyond the physical solution range to capture the diffuse tails.
The active kinetic levels are selected using
\cref{eq:active_kinetic_set} with
\(\kappa_{\mathrm{act}}=2\times10^{-3}\). This relative threshold removes
the small numerical background introduced by the diffuse lifting while
retaining the kinetic levels carrying the resolved defect.  \(\Omega_\eta\) and \(N_\eta\) are chosen to
contain and resolve the localized transverse defect profile. The resulting
domains and discretizations are summarized in \cref{tab:domain}.

\begin{table}[!htbp]
    \centering
    \begin{tabular}{llccccccc}
        \toprule
        Sec. & $\Omega_{\bm x}$&$N_{\bm x}$ & $\Omega_\xi$ &$N_\xi$& $\Omega_\eta$&$N_\eta$&$\Omega_\theta$&$N_\theta$ \\
        \midrule
\ref{subsec:ramp_riemann}&$[-2,2]$&$1000$&$[-0.1,2.1]$&$221$&$[-0.02,0.02]$&$81$& - & -\\
\ref{subsec:triangular_initial_condition}&$[-5,7]$&$1200$&$[-0.1,2.1]$&$221$&$[-0.03,0.03]$&$31$& - & -\\
\ref{subsec:burgers_2d_gaussian}&$[0,1]^2$&$1000^2$&$[-0.05,1.05]$&$111$&$[-0.12,0.12]$&$128$& $[0,1]$ & $128$\\
\ref{subsec:BL_2d}&$[0,1]^2$&$1000^2$&$[-0.05,1.05]$&$111$&$[-0.01,0.01]$&$201$& $[0,1]$ & $1000$\\
        \bottomrule
    \end{tabular}
    \caption{Domains and discretizations for the numerical examples.}
    \label{tab:domain}
\end{table}

We evaluate the ROMs primarily through the relative L$^2$ solution error,
\begin{equation}\label{eq:rel-err}
    e(t^n)
    =
    \frac{
        \left\|u_{\mathrm{ref}}(\cdot,t^n)
        -u_{\mathrm{ROM}}(\cdot,t^n)\right\|_{L^2(\Omega_{\bm x})}
    }{
        \left\|u_{\mathrm{ref}}(\cdot,t^n)\right\|_{L^2(\Omega_{\bm x})}
    }.
\end{equation}
In addition to state accuracy, we monitor two integral quantities that characterize the global conservative and entropy-dissipative behavior of the reconstructed solution: the total mass \(M(t)=\int_{\Omega_{\bm x}} u(\bm x,t)\,\mathrm d\bm x\) and the quadratic entropy \(\mathcal T(t)=\frac12\int_{\Omega_{\bm x}} u(\bm x,t)^2\,\mathrm d\bm x\). Together, these diagnostics assess whether accurate prediction of the registered defect dynamics translates into an accurate physical solution while retaining the global balance and entropy-dissipation behavior of the reference solution.

\subsection{Ramp--Riemann initial condition}
\label{subsec:ramp_riemann}
The first example considers the one-dimensional Burgers equation
\[
    \partial_t u+\partial_x\left(\frac{u^2}{2}\right)=0,\quad     u_0(x)
    =
    \begin{cases}
        2, & x\leq -1,\\
        -2x, & -1<x<0,\\
        0, & x\geq 0.
    \end{cases}
\]
The initial compression wave breaks at \(t_{\mathrm b}=0.5\), after which the entropy solution contains a single moving shock. Because the solution is available analytically, this example provides a controlled setting in which to assess whether the ROM captures both the shock trajectory and the associated entropy-producing dynamics. The analytic mass and entropy are given by 

\[ M(t)=3+2t, \qquad \mathcal T(t)=
\begin{cases}
\frac{8}{3}(1+t),&t<\frac12,\\
3+2t,&t\ge\frac12
\end{cases}\]  

The time interval is \(0\leq t\leq1.25\), with \(\Delta t=0.01\). Event snapshots begin at the breaking time, and the ROM is trained through \(t\approx0.85\); the remaining interval \(0.85<t\leq1.25\) is reserved for prediction. Since the spatial shock set consists of a single point, the shock geometry reduces to its scalar location, which is modeled by the polynomial fit described in \cref{rmk:1d}, while the registered defect-driven source is evolved by the DMD model on the \((\eta,\xi)\)-grid. A quadratic polynomial is used for the shock-location trajectory. For the registered defect-driven source, the ROM reduces the state dimension from $N_\eta N_\xi=81\times221=17{,}901$ to $r_G=2$.

\begin{figure}[!htbp]
\centering
\includegraphics[width= \textwidth]{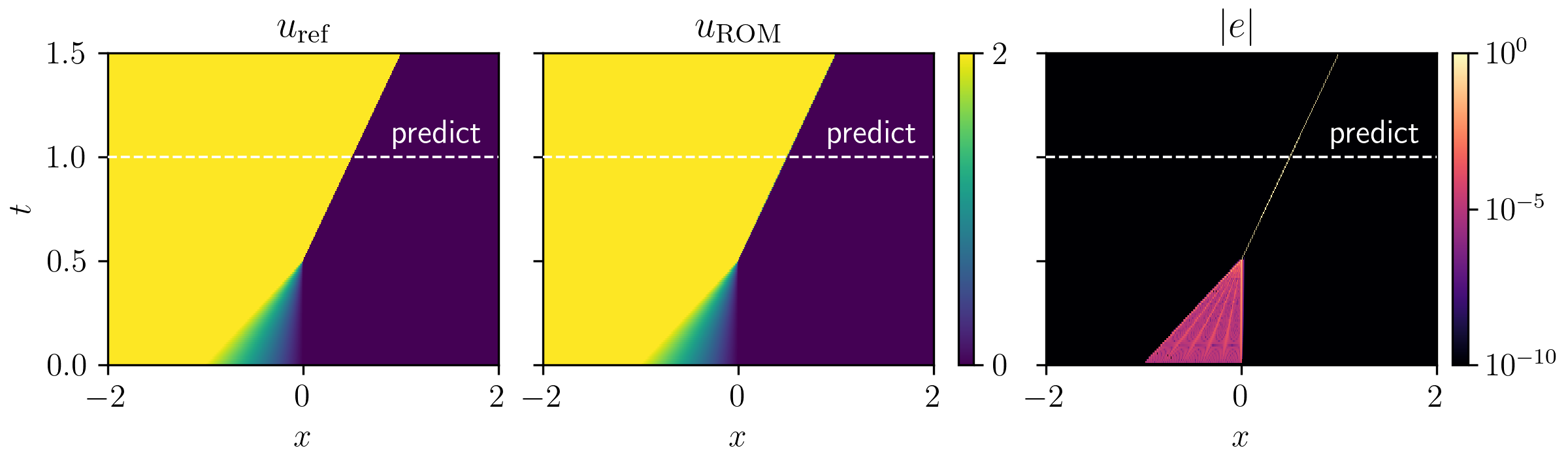}
\caption{Reference solution, ROM reconstruction, and pointwise absolute error
for \cref{subsec:ramp_riemann}. The dashed line marks the end of the training
interval.}
\label{fig:1d-burgers-riemann-solution}
\end{figure}

\begin{figure}[!h]
\centering
\includegraphics[width= \textwidth]{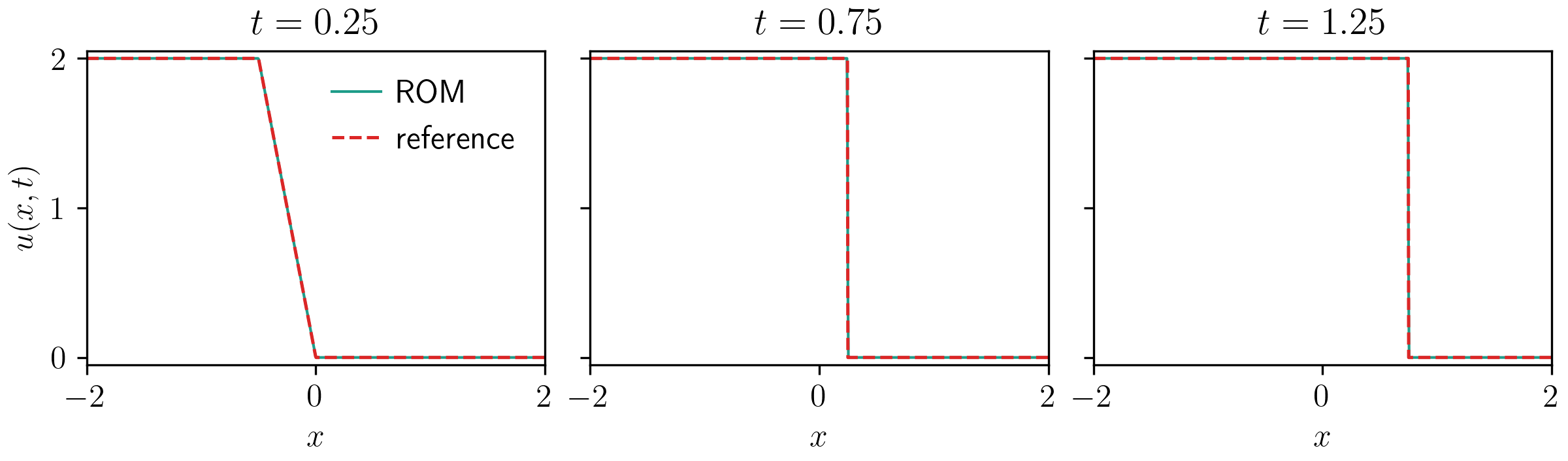}
\caption{Reference and ROM solution profiles for \cref{subsec:ramp_riemann} at selected times before and after shock formation}
\label{fig:1d-burgers-riemann-profile}
\end{figure}

\begin{figure}[!h]
\centering
\includegraphics[width= \textwidth]{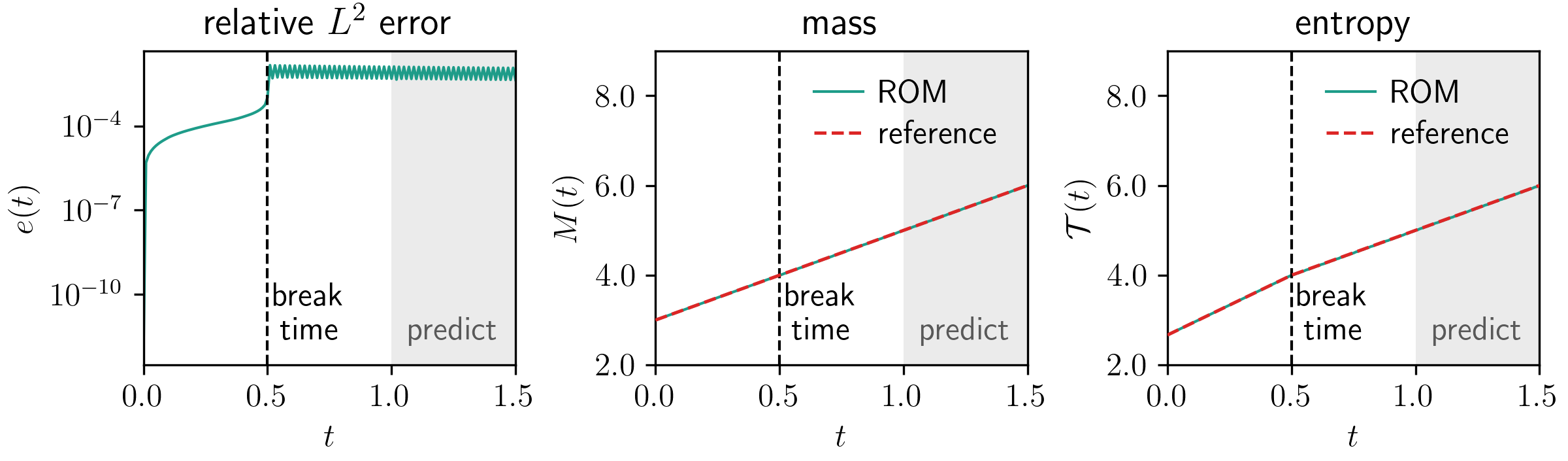}
\caption{Relative $L^2$ error, total mass, and quadratic entropy for  \cref{subsec:ramp_riemann}. The vertical dashed line marks the breaking time $t_b$,
and the prediction interval is shown in gray.}
\label{fig:1d-burgers-riemann-error}
\end{figure}

\cref{fig:1d-burgers-riemann-solution} compares the reference and \ac{ROM} space-time solutions. The ROM reproduces the pre-breaking characteristic evolution and accurately tracks the shock after $t_b$, including beyond the end of the training interval. The largest errors are concentrated near the shock, as expected for a small displacement error in a discontinuous solution. The solution profiles in \cref{fig:1d-burgers-riemann-profile} further show that the sharp shock front is accurately captured without introducing spurious oscillations.

The corresponding error and integral diagnostics are shown in \cref{fig:1d-burgers-riemann-error}. The relative L$^2$ error remains small throughout the prediction interval, while the mass and quadratic entropy closely follow the analytic reference, demonstrating accurate global balance and entropy-dissipation behavior.

\subsection{Triangular initial condition}
\label{subsec:triangular_initial_condition}
The second example uses the same one-dimensional Burgers flux with the triangular initial condition
\[
    u_0(x)
    =
    \begin{cases}
        x+3, & -3\leq x\leq-1,\\
        1-x, & -1<x\leq1,\\
        0, & \text{otherwise}.
    \end{cases}
\]
The two smooth branches meet at \(x=-1\), and the solution develops a shock at \(t_{\mathrm b}=1\). The known entropy solution is again used as the
reference with mass and entropy  given by
\[
M(t)=4,\qquad
\mathcal T(t)=
\begin{cases}
\frac83,&t<1,\\
\frac{8\sqrt2}{3\sqrt{1+t}},&t\ge1.
\end{cases}\]

Snapshots are taken with \(\Delta t=0.02\) on \(0\leq t\leq3\). The ROM is trained from \(t_b=1\) through \(t=2\) and subsequently used to predict the solution on \(2<t\leq3\). As in the first example, the shock geometry reduces to its scalar location, which is modeled by a cubic polynomial fit, while the registered defect-driven source is evolved by the DMD model, reducing the state dimension from \(N_\eta N_\xi=31\times221=6{,}851\) to \(r_G=6\).

\begin{figure}[!h]
\centering
\includegraphics[width= \textwidth]{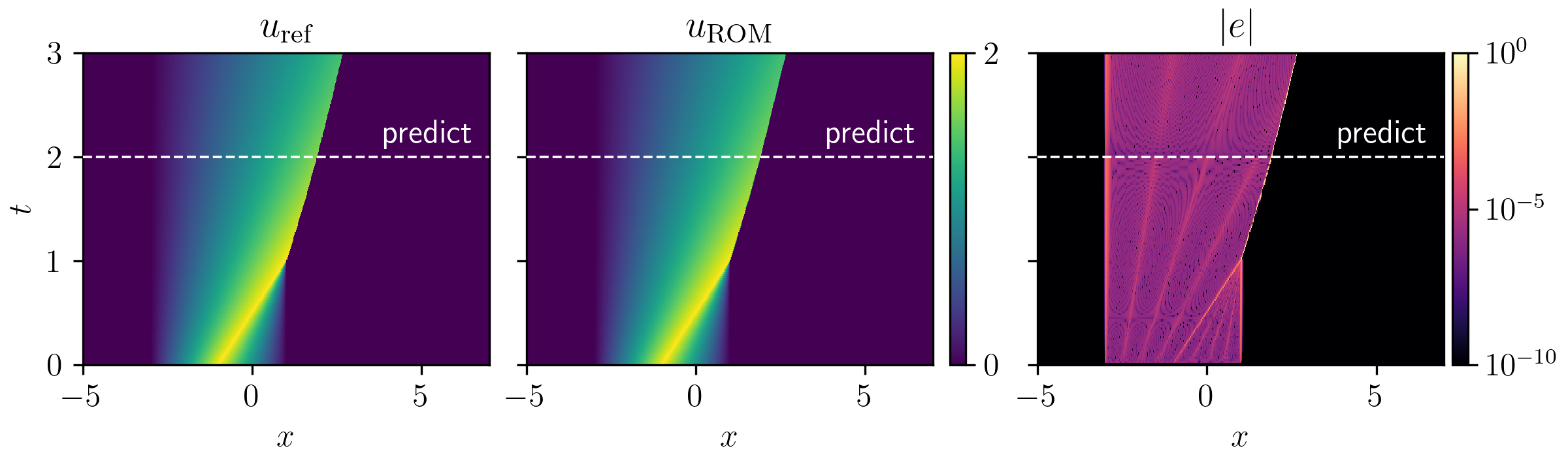}
\caption{Reference solution, ROM reconstruction, and pointwise absolute error
for \cref{subsec:triangular_initial_condition}. The dashed line marks the end of the training
interval.}
\label{fig:1d-burgers-triangle-solution}
\end{figure}

\begin{figure}[!h]
\centering
\includegraphics[width= \textwidth]{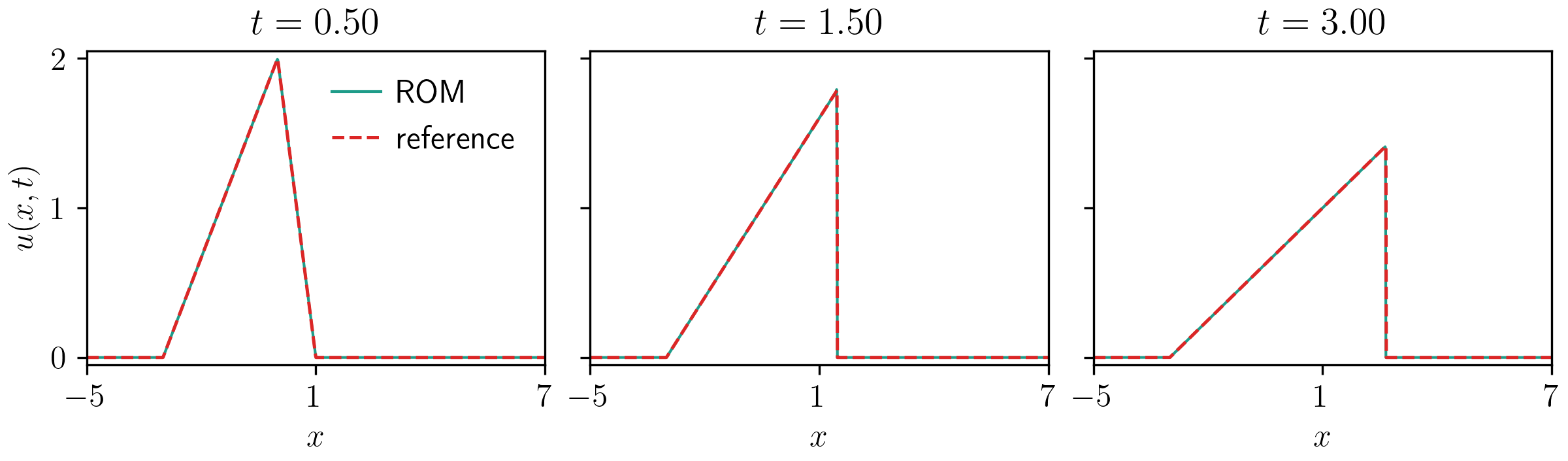}
\caption{Reference and ROM solution profiles for \cref{subsec:triangular_initial_condition} at selected times before and after shock formation}
\label{fig:1d-burgers-triangle-profile}
\end{figure}

\begin{figure}[!h]
\centering
\includegraphics[width= \textwidth]{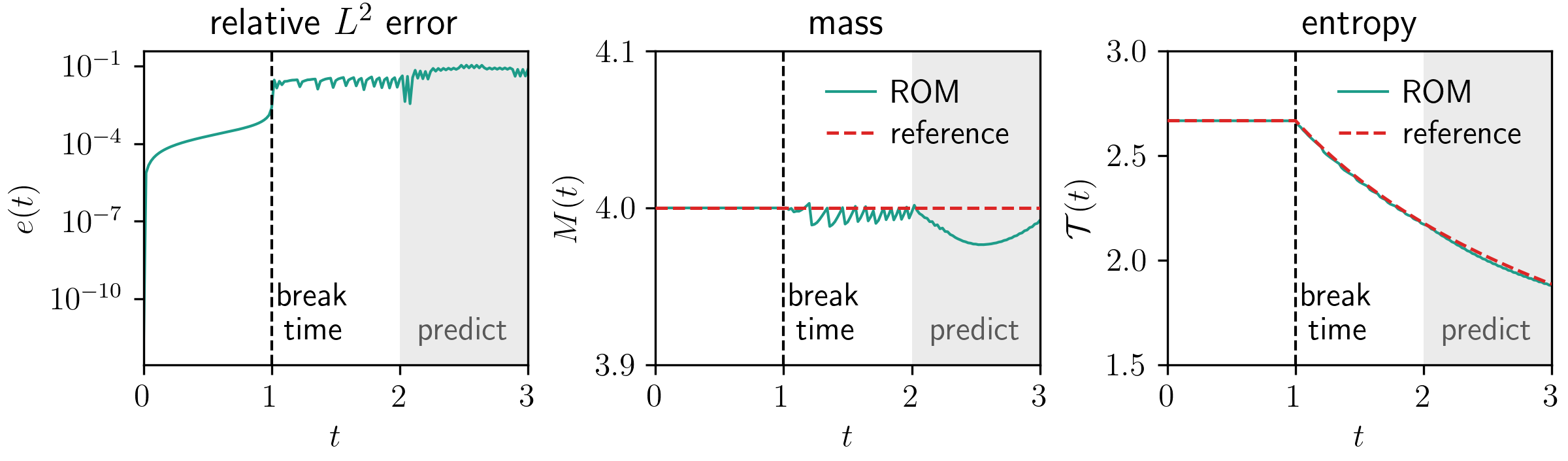}
\caption{Relative $L^2$ error, total mass, and quadratic entropy for  \cref{subsec:triangular_initial_condition}. The vertical dashed line marks the breaking time $t_b$,
and the prediction interval is shown in gray.}
\label{fig:1d-burgers-triangle-error}
\end{figure}

\Cref{fig:1d-burgers-triangle-solution} compares the reference and ROM space-time solutions. The ROM accurately captures the shock formation and subsequent propagation, with good agreement maintained beyond the end of the training interval. The largest errors remain localized near the shock front, where small errors in its predicted position produce a localized discrepancy in the discontinuous solution. The solution profiles in \cref{fig:1d-burgers-triangle-profile} further show that the ROM accurately predicts both the evolving rarefaction and the shock with time-varying propagation speed, while maintaining a sharp and smoothly connected solution without spurious oscillations. The relative L$^2$ error, mass, and quadratic entropy in \cref{fig:1d-burgers-triangle-error} show similarly good agreement throughout the prediction interval.

\subsection{Two-dimensional Burgers equation with Gaussian initial data}
\label{subsec:burgers_2d_gaussian}
The third example is a two-dimensional Burgers equation
\[
    \partial_t u
    +
    \partial_{x_1}\left(\frac{u^2}{2}\right)
    +
    \partial_{x_2}\left(\frac{u^2}{2}\right)
    =0
\]
on \(\Omega_{\bm x}=[0,1]^2\), with
\(
    u_0(x_1,x_2)
    =
    \exp\left(
        -\frac{(x_1-0.35)^2+(x_2-0.35)^2}{0.02}
    \right).
\)
Homogeneous inflow data are imposed at \(x_1=0\) and \(x_2=0\), while the remaining boundaries are outflow boundaries. The reference solution is computed by a second-order TVD Godunov discretization with minmod reconstruction and third-order SSP Runge--Kutta time integration. Snapshots are stored with \(\Delta t=0.005\) on \(0\leq t\leq0.36\). The breaking time is \(t_{\mathrm b}\approx0.12\); post-breaking snapshots from \(0.12\leq t\leq0.24\) are used for training, and the interval \(0.24<t\leq0.36\) is reserved for prediction.

The empirical defect is evaluated directly at the required physical and registered coordinates, without constructing a separate multidimensional characteristic grid. We choose 
\(\bm x_o=\bm 0,
\bm T=\frac{1}{\sqrt{2}}
\begin{bmatrix}
1\\-1
\end{bmatrix},
\bm d=\frac{1}{\sqrt{2}}
\begin{bmatrix}
1\\1
\end{bmatrix}.\) Thus, consistently with \cref{eq:local_shock_coordinates}, 
\(
    q=\frac{\sqrt{2}}{2}(x_1-x_2),
    p=\frac{\sqrt{2}}{2}(x_1+x_2)
\), with the transverse direction aligned with the characteristic velocity \(\bm f'(\xi)=(\xi,\xi)^T\). The probe domain is
\(
    \Omega_q^{\mathrm{probe}}=[-0.24,0.24],
    \Omega_p=[0.45,1.40]
\),
using \(144\) tangential and \(224\) transverse probe points. The resolved tangential support (\cref{eq:resolved_tangential_support}) is selected with \(\kappa_\Gamma=2\times10^{-2}\). The transverse maximizer of the aggregate event density is refined by local quadratic interpolation, and the resulting ridge is represented by a quadratic polynomial fit. Independent DMD models are then constructed for the shock geometry and registered defect-driven source, with their ranks \(r_c\) and \(r_G\) selected separately by the \(99\%\) snapshot-energy criterion. The corresponding state dimensions are reduced from \(N_\theta=256\) to \(r_c=1\) and from \(N_\xi N_\theta N_\eta=111\times128\times128=1{,}818{,}624\) to \(r_G=21\), respectively.

\begin{figure}[!htbp]
\centering
\includegraphics[width= \textwidth]{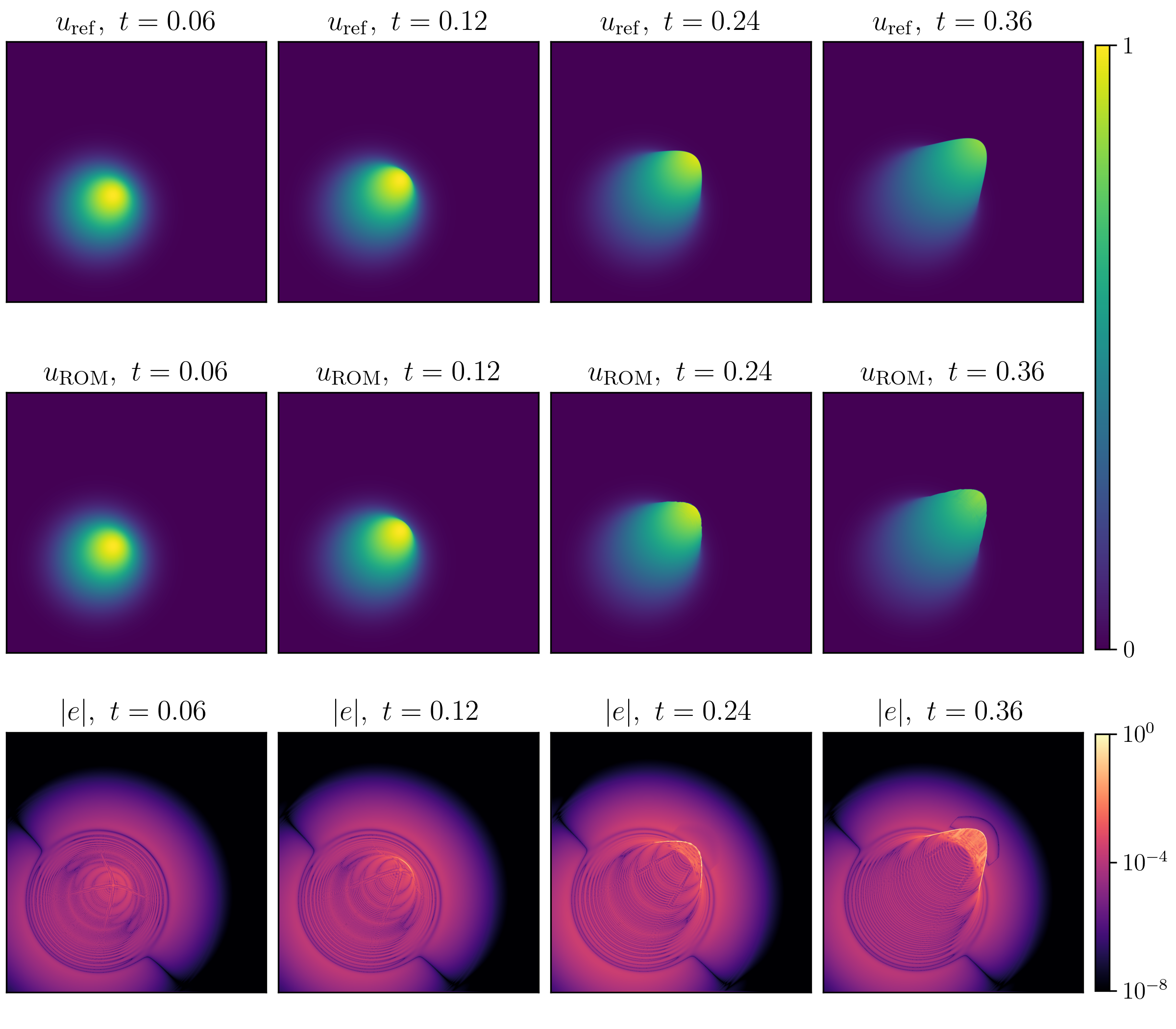}
\caption{Reference solution, ROM reconstruction, and pointwise absolute error
for \cref{subsec:burgers_2d_gaussian} at selected times before and after
shock formation.}
\label{fig:2d-burgers-solution}
\end{figure}

\begin{figure}[!htbp]
\centering
\includegraphics[width= \textwidth]{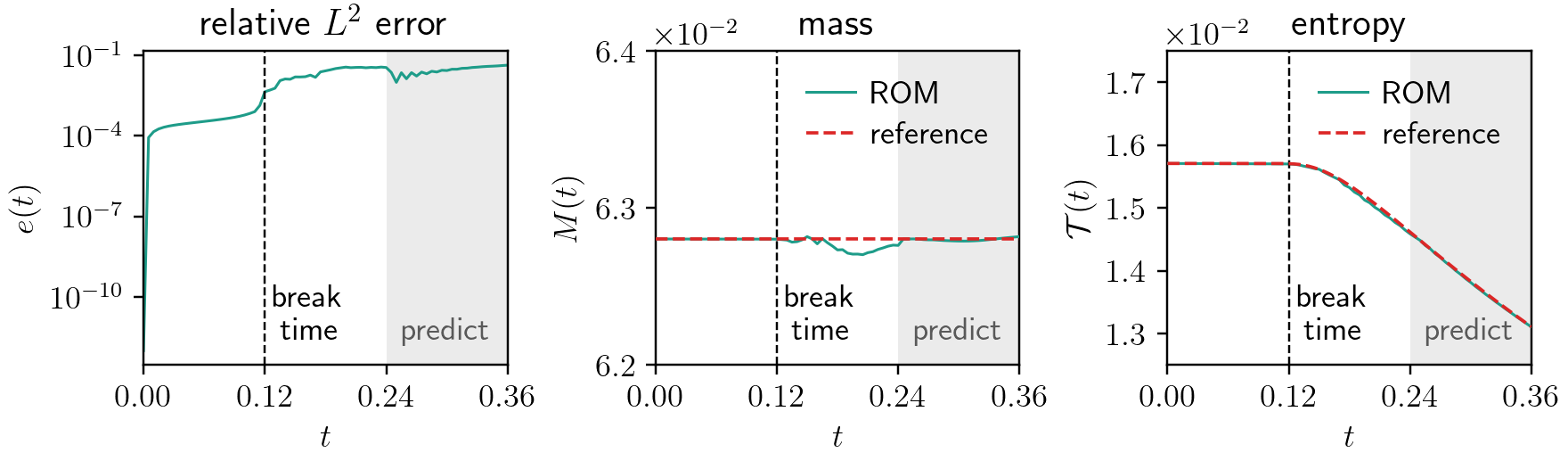}
\caption{Relative $L^2$ error, total mass, and quadratic entropy for \cref{subsec:burgers_2d_gaussian}. The vertical dashed line marks the breaking
time, and the prediction interval is shown in gray.}
\label{fig:2d-burgers-error}
\end{figure}

\Cref{fig:2d-burgers-solution} compares the reference and \ac{ROM} solutions at selected times. The \ac{ROM} accurately captures the formation and subsequent propagation of the curved shock front, with the error remaining localized near the front. The agreement is maintained beyond the training interval, demonstrating accurate prediction of both the evolving shock geometry and the associated defect-driven dynamics. The diagnostics in \cref{fig:2d-burgers-error} show that the prediction error remains controlled, while the mass and quadratic entropy closely follow the reference solution beyond the training interval.

\subsection{Two-dimensional Buckley--Leverett equation}
\label{subsec:BL_2d}

The fourth example considers two-dimensional Buckley--Leverett displacement
in a layered porous medium. The saturation satisfies
\begin{equation}
    \partial_t u
    +
    \nabla_{\bm x}\cdot
    \left(
        \bm v(\bm x)f(u)
    \right)
    =0,
    \qquad
    f(u)
    =
    \frac{u^2}{u^2+M(1-u)^2},
    \qquad
    M=2.
    \label{eq:layered_bl_problem}
\end{equation}
The porosity is set to one. We prescribe
\begin{equation}
    K(x_2)
    =
    1+\alpha\cos(2\pi x_2),
    \qquad
    \alpha=0.3,
    \qquad
    p(\bm x)=1-x_1.
    \label{eq:layered_permeability}
\end{equation}
The corresponding Darcy velocity is
\begin{equation}
    \bm v(\bm x)
    =
    -K(x_2)\nabla_{\bm x}p(\bm x)
    =
    \begin{bmatrix}
        K(x_2)\\
        0
    \end{bmatrix}.
    \label{eq:layered_velocity}
\end{equation}
Thus, the flow is directed along the \(x_1\)-axis, while its magnitude varies
with \(x_2\). The permeability satisfies
\(0.7\leq K(x_2)\leq1.3\).

The initial condition is the planar discontinuity
\begin{equation}
    u(\bm x,0)
    =
    \begin{cases}
        1, & x_1\leq0.1,\\
        0, & x_1>0.1.
    \end{cases}
    \label{eq:layered_bl_initial_condition}
\end{equation}
Unit saturation is prescribed at the left inflow boundary. Constant
extrapolation is used at the right outflow boundary and along the transverse
boundaries. Because the transport speed depends on \(x_2\), the initially
planar shock develops a smooth two-dimensional profile.

The reference solution is computed on a \(1000\times1000\) Cartesian grid
using a monotone upwind finite-volume discretization and third-order
SSP Runge--Kutta time integration. CFL-limited internal substeps are used
between stored snapshots. Snapshots are stored with
\(\Delta t=1.25\times10^{-3}\) on \(0\leq t\leq0.4\). The ROM is trained on
\(0\leq t\leq0.2\), while \(0.2<t\leq0.4\) is reserved for prediction.

For each kinetic level \(\xi\), the prescribed transport velocity is
\begin{equation}
    f'(\xi)\bm v(\bm x)
    =
    K(x_2)
    \frac{2M\xi(1-\xi)}
    {\left[\xi^2+M(1-\xi)^2\right]^2}
    \bm e_1.
    \label{eq:layered_bl_transport_velocity}
\end{equation}
The shock is represented using the graph coordinates
\begin{equation}
    \bm x_o=\bm 0,
    \quad
    \bm T=
    \begin{bmatrix}
        0\\
        1
    \end{bmatrix},
    \quad
    \bm d=
    \begin{bmatrix}
        1\\
        0
    \end{bmatrix},
    \quad
    q=\bm T^T\bm x=x_2,
    \quad
    p=\bm d^T\bm x=x_1.
    \label{eq:layered_bl_graph_coordinates}
\end{equation}
Because the shock is naturally represented as a graph over the physical \(x_2\)-coordinate, no separate graph-coordinate grid is required. We therefore set
\begin{equation}
    \theta=q = x_2,
    \qquad
    \bm c(\theta,t)
    =
    \begin{bmatrix}
         p_\Gamma(\theta,t)\\
        \theta
    \end{bmatrix},
    \label{eq:layered_bl_embedding}
\end{equation}
with \(\theta\) sharing the same discretization as \(x_2\).

\begin{figure}[!htbp]
    \centering
    \includegraphics[width=\textwidth]
        {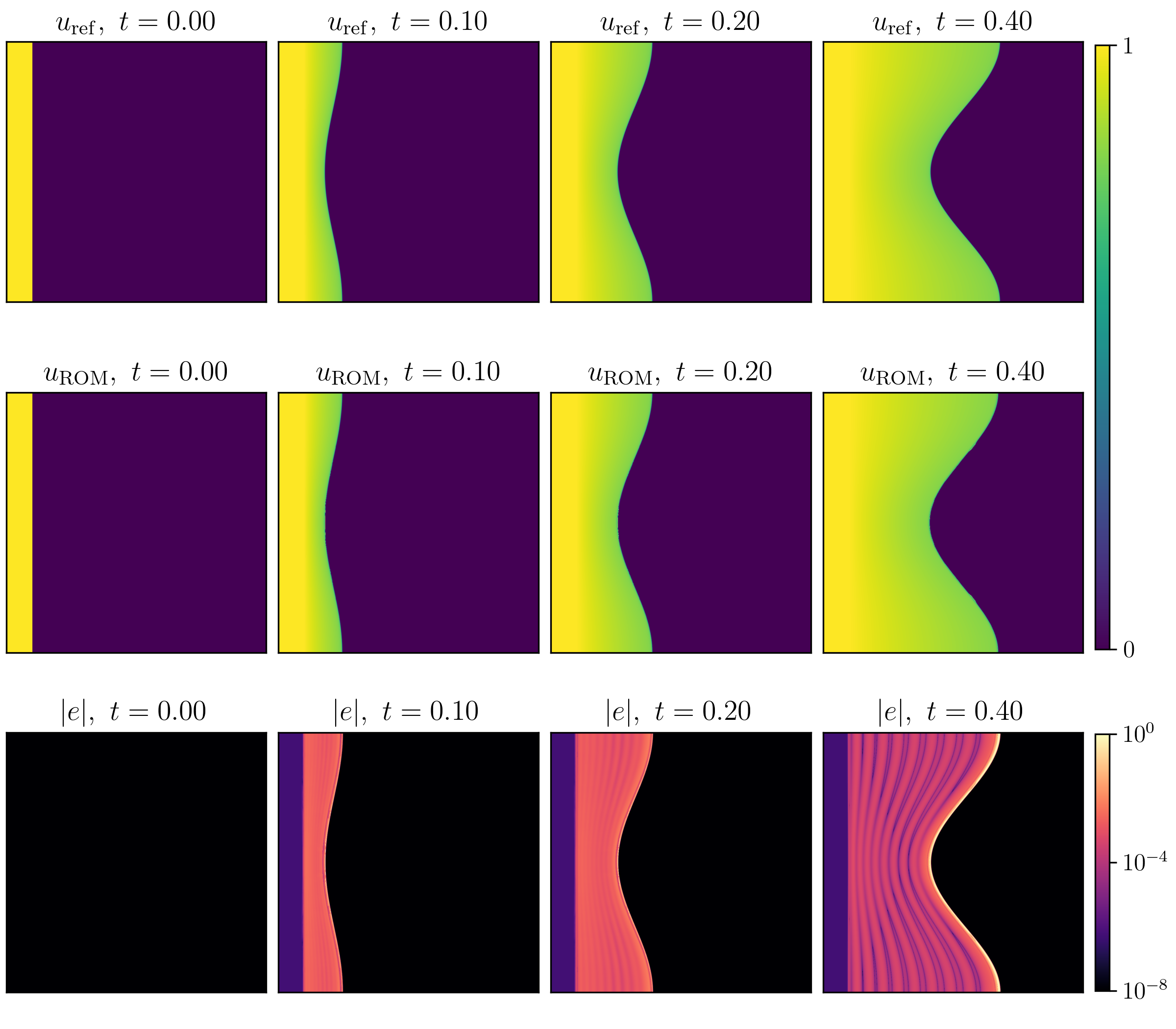}
    \caption{Reference solution, ROM approximation, and pointwise absolute
    error for \cref{subsec:BL_2d} at selected training and prediction times.}
    \label{fig:2d_bl_solution}
\end{figure}

\begin{figure}[!htbp]
    \centering
    \includegraphics[width=\textwidth]
        {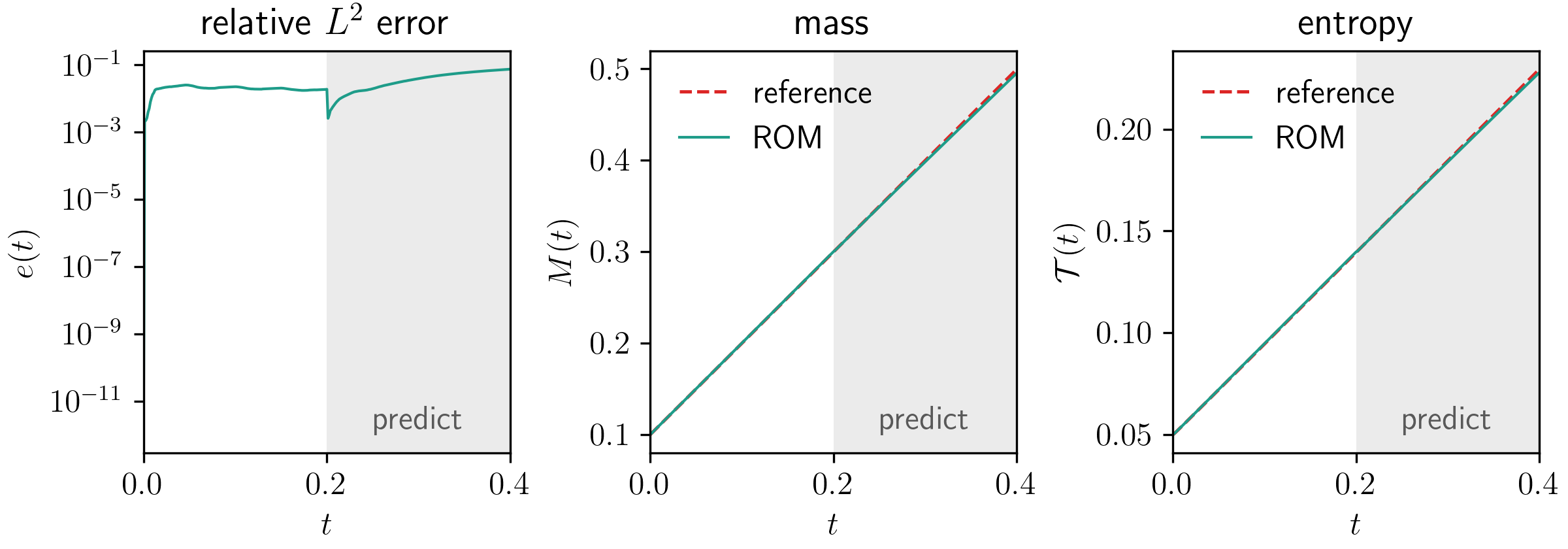}
    \caption{Relative \(L^2\) error, total mass, and quadratic entropy for
    \cref{subsec:BL_2d}. The prediction interval is shown in gray.}
    \label{fig:2d_bl_error}
\end{figure}

Following \cref{alg:kd_rom}, independent DMD models are constructed for the shock geometry and the registered defect-driven source, with their ranks selected using the \(99\%\) snapshot-energy criterion. The corresponding state dimensions are reduced from \(N_\theta=1000\) to \(r_c=2\) and from \(N_\xi N_\theta N_\eta=111\times1000\times201=22{,}311{,}000\) to \(r_G=43\), respectively.

 \cref{fig:2d_bl_solution} shows that the ROM accurately reproduces the transverse deformation and subsequent propagation of the saturation front. The largest discrepancies remain concentrated near the evolving front, where small errors in its predicted position produce localized differences in the discontinuous solution. The agreement remains strong throughout the prediction interval, demonstrating that the ROM captures both the evolving shock geometry and the associated defect-driven dynamics.

The diagnostics in \cref{fig:2d_bl_error} provide a consistent picture. The relative \(L^2\) error remains controlled beyond the training interval, while the total mass and quadratic entropy of the ROM remain in close agreement with those of the reference solution. These results demonstrate that the reduced model preserves the principal global balance properties of the full-order solution while predicting the strongly heterogeneous two-dimensional front evolution.

\section{Conclusions}
\label{sec:conclusion}
We developed a structure-informed data-driven reduced-order modeling framework for scalar hyperbolic conservation laws based on the kinetic defect formulation. The kinetic representation separates the nonlinear solution dynamics into analytically known characteristic transport and a defect-driven correction associated with entropy-producing shocks. Exploiting the localization of the kinetic defect, the proposed method extracts the evolving shock manifold and registers the empirical defect-driven source in shock-attached coordinates. Separate \acp{ROM} are then constructed for the shock geometry and the registered defect dynamics, while the known characteristic transport is prescribed directly. The numerical examples demonstrate that this decomposition enables accurate reconstruction and prediction of shock-dominated solutions beyond the training interval while retaining the global mass and entropy-dissipation behavior of the reference solutions. More broadly, these results highlight the value of incorporating physical knowledge and mathematical structure into data-driven model reduction. For challenging hyperbolic problems, such structure can guide the choice of learning variables and reduced representations, allowing the data-driven model to focus on the dynamically essential components rather than the full solution evolution. The kinetic defect formulation provides one example of how this principle can lead to effective reduced models for problems that remain difficult for conventional linear-subspace approaches.

The present work focuses on solutions whose active shock manifold can be represented within a single graph chart. An important next step is to extend the framework to problems involving multiple shocks, rarefaction waves, and wave interactions, including shock merging, splitting, and changes in front topology. Such problems will generally require multiple local graph charts together with mechanisms for identifying, evolving, and coupling the corresponding localized defect contributions as interactions occur. Developing these multi-chart representations, together with adaptive strategies for the creation and removal of charts, provides a natural direction for extending the kinetic-defect ROM framework to more complex multidimensional hyperbolic dynamics.

\bibliographystyle{elsarticle-num}
\bibliography{main}
\end{document}